\documentclass[11pt]{article}
\usepackage{amssymb,latexsym,amsmath}
\usepackage{amsfonts,amsbsy,bm}

\def\a {\alpha }

\newcommand {\SB} {{\Bbb B}}
\newcommand {\SC} {{\Bbb C}}

\newcommand {\SN} {{\Bbb N}}
\newcommand {\SR} {{\Bbb R}}

\newcommand {\inRn} {\in\SR^n}

\renewcommand {\phi} {{\varphi}}

\newcommand {\al} {{\alpha}}

\newcommand {\Dt} {{\Delta}}
\newcommand {\e} {{\varepsilon}}
\newcommand {\ga} {{\gamma}}

\newcommand {\la} {{\lambda}}

\renewcommand {\O} {{\Omega}}

\newcommand {\cB} {{{\cal B}}}

\newcommand {\cD}   {{\Tu}}  
\newcommand {\cE} {{\cal E}}
\newcommand {\cF} {{\cal F}}

\newcommand {\cH} {{\cal H}}

\newcommand{\cN}{{\cal N}}

\newcommand {\cS} {{\cal S}}

\newcommand {\tF} {{\widetilde F}}

\newcommand {\tpn} {{{\tilde{p}_\nu}}}

\newcommand {\hf} {{\widehat f}}

\newcommand {\hpsi} {{\widehat\psi}}

\newcommand {\hphi} {{\widehat\phi}}

\renewcommand {\tfrac}[2]{{\textstyle\frac{#1}{#2}}}

\newcommand {\supp} {\mbox{Supp }}

\newcommand {\Scirc} {\raise.2ex\hbox{$\scriptstyle\circ$}}

\newcommand {\im} {{\Im\!\mbox{\small\it m}\,}}

\newcommand {\mand} {{\quad\mbox{and}\quad}}
\renewcommand {\mid} {{\,\,\,\colon\,\,\,}}
\newcommand {\Ol} {\overline}

\newcommand {\Proof} {\noindent{\bf P{\footnotesize\bf ROOF}: } \ }
\newcommand {\Proofof}[1] {\noindent{\bf P{\footnotesize\bf ROOF} of {#1}: } \ }
\newcommand {\ProofEnd} {
             \begin{flushright} \vskip -0.2in $\Box$ \end{flushright}}

\newcommand{\Ba}[1]{\begin{array}{#1}}
\newcommand{\Ea}{\end{array}}
\newcommand{\Be}{\begin{equation}}
\newcommand{\Ee}{\end{equation}}
\newcommand{\Bea}{\begin{eqnarray}}
\newcommand{\Eea}{\end{eqnarray}}
\newcommand{\Beas}{\begin{eqnarray*}}
\newcommand{\Eeas}{\end{eqnarray*}}
\newcommand{\Benu}{\begin{enumerate}}
\newcommand{\Eenu}{\end{enumerate}}
\newcommand{\Bi}{\begin{itemize}}
\newcommand{\Ei}{\end{itemize}}

\newcommand{\BR}{\begin{Remark} \em}
\newcommand{\ER}{\end{Remark}}
\newcommand{\BE}{\begin{example} \em}
\newcommand{\EE}{\end{example}}
\newcounter{remark}

\newtheorem{theorem}[equation]{T{\hskip 0pt\footnotesize\bf HEOREM}}
\newtheorem{proposition}[equation]{P{\hskip 0pt\footnotesize\bf ROPOSITION}}
\newtheorem{corollary}[equation]{C{\hskip 0pt\footnotesize\bf OROLLARY}}
\newtheorem{lemma}[equation]{L{\hskip 0pt\footnotesize\bf EMMA}}
\newtheorem{Remark}[equation]{R{\hskip 0pt\footnotesize\bf EMARK}}
\newtheorem{definition}[equation]{D{\hskip 0pt\footnotesize\bf EFINITION}}
\newtheorem{example}[equation]{E{\hskip 0pt\footnotesize\bf XAMPLE}}
\newtheorem{conjecture}{C{\hskip 0pt\footnotesize\bf ONJECTURE}}

\newcommand{\bO}{\partial\Omega}

\newcommand{\Tu}{{T_\O}}

\newcommand{\be}{{\bf e}}

\newcommand{\fnr}{\frac{n}r}

\newcommand{\cLL}{{\cal L}}

\renewcommand{\O}{\Omega}

\newcommand{{\tBox}}{{\widetilde{\raisebox{-0.2ex}[1.25ex][0ex]{$\Box$}}}}

\newcounter{reb}
\newcounter{rea}
\newcounter{red}
\newcounter{reg}
\newcounter{res}
\newcommand{\bline}{{\bigskip

\noindent}}

\newcommand{\sline}{{\smallskip

\noindent}}

\begin{document}
\title{Analytic Besov spaces and Hardy-type inequalities in tube
  domains over symmetric cones.\footnotetext{\emph{2000 Math Subject Classification:} 42B35, 32M15.}
}
\author{D. B\'ekoll\'e\quad A. Bonami$^*$\quad
G. Garrig\'os \footnote{Research partially supported by the
European Commission, within the IHP Network ``HARP 2002-2006'',
contract number HPRN-CT-2001-00273-HARP. Third author also
supported by \emph{Programa Ram\'on y Cajal} and grant
``MTM2007-60952'', MEC (Spain).} \quad F. Ricci$^*$ \quad B.
Sehba} \maketitle

\begin{abstract}
We give various equivalent formulations to the (partially) open
problem about $L^p$-boundedness of Bergman projections in tubes
over cones. Namely, we show that such boundedness is equivalent to
the duality identity between Bergman spaces, $A^{p'}=(A^p)^*$, and
also to a Hardy type inequality related to  the wave operator. We
introduce analytic Besov spaces in tubes over cones, for which
such Hardy inequalities play an important role. For $p\geq 2$ we
identify  as a Besov space the range of the Bergman projection
acting on $L^p$, and also the dual of $A^{p'}$. For the Bloch
space $\SB^\infty$ we give in addition new necessary conditions on
the number of derivatives required in its definition.
\end{abstract}

\section{Introduction}
\setcounter{equation}{0} \setcounter{footnote}{0}
\setcounter{figure}{0}

 Let $\cD$ be a symmetric domain of tube
type in $\SC^n$, that is $\cD=\SR^n+i\O$ where $\O$ is an
\emph{irreducible symmetric cone} in $\SR^n$. These domains can be
seen as  multidimensional analogues of the upper half plane in
$\SC$. A typical example arises when $\O$  is the forward
light-cone of $\SR^n$, $n\geq3$,
\[
\Lambda_n\,=\,\bigl\{y\in\SR^n\mid
y_1^2-y_2^2-\ldots-y_n^2>0,\;\;y_1>0\bigr\}.
\]
Other examples correspond to the cones $\mbox{Sym}_+(r,\SR)$ of
 positive definite symmetric $r\times
r$-matrices. We refer to the text \cite{FK} for a general
description of symmetric cones. Following the notation in
\cite{FK} we write $r$ for the rank of $\O$ and $\Dt(x)$ for the
associated determinant function. In the above examples,
light-cones have rank $2$ and determinant equal to the Lorentz
form $\Dt(y)= y_1^2-y_2^2-\ldots-y_n^2$, while the cones
$\mbox{Sym}_+(r,\SR)$ have rank $r$ and the determinant is the
usual determinant of $r \times r$ matrices. We shall denote by
$\mathcal{H}(\cD)$ the space of holomorphic functions on $\cD$.

 A major open
question in these domains concerns the $L^p$ boundedness of
\emph{Bergman projections}, which can only hold for values of $p$
sufficiently close to 2 \cite{BB1,BBPR,BBGR}. More precisely,
consider the (weighted) spaces
\[L^p_\nu(\cD)\,=\,L^p(\cD,\Dt(y)^{\nu-n/r}dx\,dy)\] and let $A^p_\nu(\cD)$
be the subspace of holomorphic functions. Denote by $P_\nu$ the
orthogonal projection mapping $L^2_\nu(\cD)$ into $A^2_\nu(\cD)$.
The usual (unweighted) Bergman spaces correspond to $\nu=\frac
nr$, while the weighted cases can be considered when $\nu>\frac
nr-1$ (since otherwise $A^p_\nu=\{0\}$).
\begin{conjecture}\label{critical}
Let $\nu>\frac nr-1$. Then the Bergman projection $P_\nu$ admits a
bounded extension to $L^p_\nu(\cD)$ if and only if
\[
p'_\nu\,< \,p\,<\,p_\nu:=\,\frac{\nu+\frac {2n}r-1}{\frac
nr-1}-\,\frac{(1-\nu)_+}{\frac nr-1}.
\]
\end{conjecture}

This problem has only been settled in the case of light-cones for
sufficiently large $\nu$'s \cite{BBGR}. In general, the known
results can be described as follows (see
\cite{BB1,BBPR,BBG,BBGR}). The fact that boundedness can only hold
when $\tpn'<p<\tpn$, where
$$\tilde p_\nu:= \,\frac{\nu+\frac {2n}r-1}{\frac
nr-1},$$ is trivially given by the $L^{p'}_\nu$-integrability of
the Bergman kernel (which only happens when $p<\tpn$) and duality.
The necessity of the condition involving $(1-\nu)_+$ was
established in \cite{BBGR}, and may only occur in the three
dimensional forward light-cone (the only case in which $\nu$ is
allowed to take values below $1$). Concerning sufficiency, it has
been proved in \cite{BBPR,BBG} that $P_\nu$ is bounded in
$L^p_\nu$ at least in the range \Be\bar p_\nu'<p< \bar p_\nu:=
\,\frac{\nu+\frac {2n}r-2}{\frac nr-1}.\label{barpnu}\Ee In the
light-cone setting (ie when $r=2$), Conjecture 1 is  closely related to other deep
conjectures for the wave equation. As shown in \cite{BBGR}, this
implies slight improvements in the range \eqref{barpnu} for all
$\nu$'s, and in fact sets completely the conjecture when $\nu$ is
sufficiently large (see also \cite{GS,GSS} for the latest
results).

In this paper, we shall not improve these boundedness results, but
interest ourselves
 in equivalent
 formulations of Conjecture 1 and implications in the theory of holomorphic function spaces in $\cD$.
Consider the ``box operator'' of $\O$, denoted $\Box=\Dt(\frac 1i
\frac{\partial}{\partial x})$, as the differential operator of
degree $r$ in $\SR^n$ defined by the equality: \Be
\Box\,[e^{i(x|\xi)}]=\Dt(\xi)e^{i(x|\xi)}, \quad x,\xi\inRn.
\label{bbox} \Ee In the rank 1 setting (that is, when $n=1$ and
$\O=(0,\infty)$) this corresponds to $-i\frac{d}{dx}$, and in the
rank 2 situation (that is, when $\O$ is the forward light cone in
$\SR^n$) we have
$\Box=-(\partial^2_{x_1}-\partial^2_{x_2}-\ldots-\partial^2_{x_n})/4$,
which explains why $\Dt(\frac{\partial}{\partial x})$ is sometimes
called the wave operator. We denote by $\Box_z=\Dt(\frac 1i
\frac{\partial}{\partial z})$ the corresponding differential
operator in $\SC^n$ defined replacing $x$ in \eqref{bbox} by
$z\in\SC^n$. Observe, however, that $\Box_z=\Box_x$ when acting on
holomorphic functions in $\cD$. To simplify notation, we will
write $\Box$ instead of $\Box_z$. Our first result can then be
stated as follows.

\begin{theorem} \label{hardyth}
Let $\nu>\frac nr-1$. Then, for $p\geq 2$, the Bergman projection
$P_\nu$ admits a bounded extension to $L^p_\nu(\cD)$ if and only
if there exists a constant $C$ such that, for all $F \in A^p_\nu$
we have
\begin{equation}\label{hardy}
\int\!\!\int_\cD |F(x+iy)|^p \,\Delta^{\nu - {n \over r}}
(y)\,dx\,dy \, \leq\, C\,\int\!\!\int_\cD \,\bigl|\Delta (y) \Box
 F(x+iy)\bigr|^p\,\Delta^{\nu - {n \over r}} (y)\,dx\,dy.
\end{equation}
\end{theorem}
We will refer to \eqref{hardy} as  \emph{Hardy inequality (for the
parameters $(p, \nu)$)}, by reference to the one dimensional
setting $n=r=1$, where it is true for all $\nu>0$ and $1\leq
p<\infty$. More comments on Hardy inequalities
 for holomorphic functions in $\cD$ have been done in \cite{B}, where a weaker statement
 was announced (see also \cite{BBPR}).


 We remark that \eqref{hardy} is always valid when
$1\leq p\leq 2$, as can be proved, for instance, from an explicit
formula for $F$ in terms of $\Box F$ involving the fundamental
solution of the Box operator (see \cite{B}). However, in this
range \eqref{hardy} has no implications in terms of boundedness of
Bergman projections. We also remark that the converse inequality,
\begin{equation}\label{mean}
\bigl\|\Dt(\im\cdot)\Box F\bigr\|_{L^p_\nu}\,\leq\,C\,\bigl\|
F\bigr\|_{L^p_\nu}
\end{equation}
for $F\in\cH(\cD)$, is valid for all $0<p\leq\infty$ and
$\nu\in\SR$, and is an easy consequence of the mean value
inequality for holomorphic functions (see \cite{BBPR}).  We will
prove Theorem \ref{hardyth} in Section 3, and add more comments on
Hardy inequalities.
\medskip

The second equivalent formulation of Conjecture 1 concerns
duality.
\begin{theorem}\label{duality}
Let $\nu > \frac nr-1$ and $1<p < \infty.$ Then $P_\nu$ admits a
bounded extension to $L^p_\nu(\cD)$ if and only if the natural
mapping of $A^{p'}_\nu$ into $(A^p_\nu)^*$  is an isomorphism.
\end{theorem}

We prove a bit more, if $p>\tpn'$ then the inclusion $\Phi:A^{p'}_\nu\hookrightarrow (A^p_\nu)^*$ is
injective, and hence boundedness of $P_\nu$ is actually equivalent
to surjectivity of $\Phi$. When $p\geq\tpn$ these two properties
fail, and $(A^p_\nu)^*$ is a space strictly larger than
$A^{p'}_\nu$ which we do not know how to identify. When $1\leq
p\leq 2$, however, it is always possible to identify $(A^p_\nu)^*$
as a ``Besov space'' of analytic functions modulo equivalence
classes, which we do in section 4. Equivalence classes appear
naturally in this setting since the injectivity of $\Phi$ (or
equivalently of $\Box |_{A^{p'}_\nu}$) fails when $p<\tpn'$. We do
not know whether in this range $\Phi$ or $\Box$ may be surjective,
a question not considered before to which we will come back later.

\medskip

In section 4 we develop the theory of analytic
Besov spaces. These arise naturally in an attempt to give a
meaning to $(A^p_\nu)^*$ or $P_\nu(L^p_\nu)$ for indices $p, \nu$
for which the operator $P_\nu$ is unbounded (see eg the one
dimensional theory in \cite{Z1}). In addition, their definition is
very closely linked with the validity of Hardy inequalities, and
for this reason we take up this matter here, leaving to subsequent
works the development of further properties.
 It is remarkable that one can develop most of this theory
without making use of the (conceptually more complicated)
\emph{real variable} Besov spaces adapted to the cone, which were
introduced in \cite{BBGR}.

To be more precise, for $\nu\in\mathbb{R}$ and $1\leq
p\leq\infty$, we define
\Be\mathbb{B}^p_\nu(\cD):=\{F \in\mathcal{H}(\cD)\mid \Dt^k(\im
.)\Box^k F\in L^p_\nu \}\label{bem}\Ee for a large enough integer
$k \geq k_0(p,\nu)$ to be given later. This definition is similar
to the one
dimensional setting \cite{Flett}
, with the role of complex derivative now played by the operator
$\Box$. The best choice of the value $k_0(p,\nu) $ is related to
the validity of Hardy inequality for $(p,\nu+pk_0)$, since only in
this case we can guarantee the equivalence of norms for different
$k$'s. Of course, when $k$ can be taken equal to $0$ one has
$\mathbb{B}^p_\nu=A^p_\nu$, but in general one must deal with
equivalence classes modulo holomorphic functions annihilated by
$\Box^k$. This is a new (and sometimes disturbing) feature
compared to the theory of analytic Besov spaces in bounded
symmetric domains developed by K. Zhu \cite{Z2}. When $p=\infty$,
the analytic Besov space $\SB^\infty$ is the usual \emph{Bloch
space} (see e.g. \cite{bek,Be}).

 Among our results we shall prove the following. Here $P^{(k)}_\nu(f)$ denotes the equivalence class $P_\nu(f)+\ker\Box^k$ (defined at least for $f$ in the dense set $L^2_\nu\cap L^p_\mu$).

\begin{theorem}\label{th3}Let $\nu>\frac
nr-1$, $2\leq p\leq\infty$ and $k\geq k_0(p,\nu)$. Then
\medskip

\noindent 1.- For every real $\mu\leq\nu$, the operator $P^{(k)}_\nu$ extends continuously  from $L^p_\mu$ onto $\SB^p_\mu$.
\medskip

\noindent 2.- The dual space
$(A^{p'}_\nu)^*$ identifies with $\SB^p_\nu$, under the pairing
\[ \langle F,G\rangle_{\nu,k}=\int_{\cD}F(z)\,\Dt^{k}(\im z)\,\overline
{\Box^{k}G(z)} \,dV_\nu(z),\quad F\in A_{\nu}^{p'},\quad G\in
\SB_{\nu}^{p}.\]
\end{theorem}

These properties are standard in the Bergman space theory of
bounded symmetric domains (see eg \cite{Z1,Z2}), as far as one
allows to take $k$ sufficiently large. The point here is to find
the smallest number of derivatives in \eqref{bem} so that these
hold. As mentioned above, this is a non trivial question directly
related with Conjecture 1.

We will be more precise about this point: if Conjecture 1 holds,
then Theorem \ref{hardyth} implies that $\SB^p_\nu$ is independent
of $k$ (and Theorem \ref{th3} is true) whenever \Be k+\tfrac\nu
p>\max\;\bigl\{(\tfrac nr-1)\tfrac1p, \;(\tfrac
nr-1)(1-\tfrac2p)-\tfrac1p, \;(\tfrac
nr-1)(\tfrac12-\tfrac1p)\;\bigr\}. \label{mnup}\Ee Thus, one can
conjecture that \eqref{mnup} defines the smallest integer for
which the above properties hold. With the presently known results
(ie the boundedness of $P_\nu$ in the range \eqref{barpnu}) we are
constrained to consider larger integers, namely numbers $k$ so
that\Be k+\tfrac\nu p>\max\;\bigl\{(\tfrac nr-1)\tfrac1p,
\;(\tfrac nr-1)(1-\tfrac2p)\;\bigr\}, \label{mnup2}\Ee which is
the same condition as \eqref{mnup} only when $1\leq p\leq 3$
(i.e., when the maximum in \eqref{mnup2} is attained at the first
number, and a bit more than this in the case of light-cones), or
when $p=\infty$. We also observe that the best integer $k$
satisfying \eqref{mnup2} is at most one unit above the optimal
integer for \eqref{mnup}.

Related with this question
 one can also consider a weaker
property  than Hardy's inequality (but apparently as difficult);
namely
\medskip

 \noindent{\bf
Question:} \emph{Given $1\leq p\leq\infty$ and $\nu\in\SR$, find
the smallest $\ell=\ell(p,\nu)\in\SN$ so that, for all $m\geq1$,
\Be \inf_{{H\in\cH(\cD)}\mid
\Box^{\ell+m}H=0}\,\bigl\|\Dt^\ell\Box^\ell(F+H)\bigr\|_{L^p_\nu}\,
\lesssim\,\bigl\|\Dt^{\ell+m}\Box^{\ell+m}F\bigr\|_{L^p_\nu},
\label{Gconj}\Ee for all holomorphic $F$ for which the right hand
side is finite.} \medskip

When $\ell=0$, this is equivalent to the surjectivity of
$\Box^m:A^p_\nu\to A^p_{\nu+mp}$, that is, whether $\Box^m F=G$
may have some solution $F\in A^p_\nu$ when the datum $G\in
A^p_{\nu+mp}$.

Hardy's inequality for $(p,\nu+\ell p)$ easily implies
\eqref{Gconj}, which hence holds in the range \eqref{mnup2} (with
$k$ replaced by $\ell)$. However, we do not know whether the
converse may be true. In fact, we do not even know whether
\eqref{mnup} is a necessary condition for \eqref{Gconj}. Below we
shall prove that the integer $\ell$ at least must satisfy
\[
\ell+\tfrac\nu p>\max\;\bigl\{(\tfrac nr-1)\tfrac1p, \;(\tfrac
nr-1)(\tfrac12-\tfrac1p)\;\bigr\}.
\]
 We remark that these type of necessary
conditions had not been considered at all in previous work. For
instance, for the Bloch space, one can ask whether there exist
functions $F\in\SB^\infty$ so that
$\|\Dt^j\Box^j\tF\|_\infty=\infty$ for all $j\leq \fnr-1$ and all
$\tF=F$ (mod Ker $\Box^{k_0}$) where $k_0=\lceil\fnr-1\rceil$; in
such case $k_0$ would really be a critical number of equivalence
classes. The classical example $F(z)=\ln(z\cdot\be+i)\in
\SB^\infty$ only has this property  for $j=0$, and it does not
seem easy to produce explicit examples with $j\geq1$. See however
Proposition \ref{NCbloch} below for the existence of such
functions with $j\leq(\fnr-1)/2$.

\medskip

Returning to the complex Besov spaces $\SB^p_\nu$, in section 4.4
we present a real variable characterization in terms of
``Littlewood-Paley decompositions'' of the cone, as described in
\cite{BBGR}. Roughly speaking, functions $F\in\SB^p_\nu$ have
Shilov boundary values $f=\lim_{{y\to 0}\atop{y\in\O}}F(x+iy)$
which are distributions in $\SR^n$, with Fourier transform
supported in $\overline\O$ and satisfying a growth condition\[
\Bigl\{\Delta^{-\frac\nu p} (\xi_j)\,||f
* \psi_j||_p\Bigr\}\in\ell^p,
\] for a suitable partition of unity $\{\psi_j\}$ associated with a lattice set $\{\xi_j\}$ of $\O$.
Conversely, every such distribution can  be extended via
Fourier-Laplace transform into a holomorphic function in
$\SB_\nu^p$. This allows in some cases to improve the value of $k$
for which the elements of the Besov space can be identified with
equivalence classes modulo holomorphic functions annihilated by
$\Box ^k$. In addition, we consider the real version of Bloch
spaces (which is new), and use this characterization to prove the
necessary conditions for \eqref{Gconj} alluded above.

\medskip


Finally, we mention the special family of Besov spaces
corresponding to the weight $\nu=-n/r$ in \eqref{bem}; that is,
$$\mathbb{B}^p\,=\,\bigl\{\,F \in\mathcal{H}(\cD)\mid \Dt^k(\im .)\Box^k F\in L^p(d\la)\,\bigr\}.$$
Here $d\lambda={\Dt^{-\frac{2n}r}(y)}dx\,dy$ denotes the invariant
measure under conformal transformations of $\cD$. These are the
analog for $\Tu$ of the Besov spaces introduced by Arazy and Yan
in bounded symmetric domains \cite{A,Y,Y1}. Special properties of
these spaces, such as M\"obius invariance and characterizations of
(small) Hankel operators will be described in subsequent papers
\cite{G,S,BGNS}.

\medskip

 The paper is structured as follows: in section 2 we
present some prerequisites about cones and Bergman kernels. In
section 3 we prove Theorems \ref{hardyth} and \ref{duality}. In
section 4 we introduce Besov and Bloch spaces and prove Theorem
\ref{th3}. The real analysis characterization is in $\S\S$4.4 and
4.5 and the necessary conditions related with \eqref{Gconj} are in
$\S$4.6. Finally, section 5 contains a brief list of open
questions which we could not answer in relation with this topic.
Besides Conjecture 1, the main problem that we leave open concerns
the question in \eqref{Gconj}.

\section{Bergman kernels and reproduction formulas}
\setcounter{equation}{0}
\subsection{Some prerequisites}

Below we shall use some invariance properties of determinants and
Box operators. To introduce them we need to recall some basic
 facts about symmetric cones (see the text
\cite{FK}).

Considering $V=\SR^n$ as a Jordan algebra, we denote its unit
element by $\be$ (think of the identity matrix in the cone of
positive definite symmetric matrices, or the point $\be=(1,
{\bf0})$ in the forward light cone). Let $G$ be the identity
component of the group of invertible linear transformations which
leave the cone $\O$ invariant. It is well known that $G$ acts
transitively on $\Omega$, which may be identified with the
Riemannian symmetric space $G/K$, where $K$ is the compact
subgroup of elements of $G$ which leave $\be$ invariant.
The determinant function is also preserved by $G$, in such a way
that \Be\label{inv-det}
\Dt(gy)=\Dt(g\be)\Dt(y)=\mbox{Det}(g)^{\frac rn}
\Dt(y),\quad\forall\;g\in G,\;y\in\O.\Ee It follows from this
formula that an invariant measure in $\O$ is given by
$\Dt(y)^{-\frac{n}r}\,dy$. The invariance of the Box operator
through the action of $G$ is an easy consequence of its definition
and the invariance of the determinant function, namely
\Be\label{inv-box} \Box
\bigl[F(g\cdot)\bigr]\,=\,\Dt(g\be)\,\bigl[\Box
F\bigr](g\cdot)\,=\,\mbox{Det}(g)^{\frac rn}\bigl[\Box
F\bigr](g\cdot),\quad\forall\;g\in G.\Ee Another fundamental
property is the following \cite[p. 125]{FK}: for every $\alpha\in
\mathbb{R} $ one has the identity in $\O$
\begin{equation}\label{boxdelta}
\Box \Delta^{\alpha}=b(\alpha)\Delta^{\alpha-1}
\end{equation}
where $b(\alpha)$ vanishes only for the $r$ values $0, \alpha_0,
\cdots (r-1)\alpha_0$, where $\alpha_0= -\frac{\frac nr -1}{r-1}$.
In particular,
\begin{equation}\label{zero}
\Box \Delta^{-\frac nr +1}(y)=0,\quad y\in\O.
\end{equation}
\subsection{Bergman kernels and Determinant function}
The (weighted) Bergman projection $P_\nu$ is defined by
$$P_\nu F(z)=\int_{\cD}B_\nu(z, w) F(w)
dV_\nu(w),
$$
where
 $B_\nu(z, w)=
 c_\nu\,\Dt^{-(\nu+\fnr)}((z-{\Ol{w}})/i)$
is the reproducing kernel of $A^2_\nu$, which we shall call
Bergman kernel (see \cite{FK}). For simplicity, we have written
$dV_\nu(w):=\Delta^{\nu-\frac nr}(v) du\,dv$, where $w=u+iv$ is an
element of $\cD$. Observe from \eqref{boxdelta} that
\Be\Box^m_z\left[B_\nu(z-\bar
w)\right]\,=\,c_{\nu,m}\,B_{\nu+m}(z-\bar w)\label{boxB} \Ee for a
suitable constant $c_{\nu,m}$, and all $m\in\SN$.
We will need integrability properties of  the determinants and
Bergman kernels, which are given by the next lemma.

\begin{lemma}\label{int-compl} Let $\a,\nu$ be real and $p>0$.
Then
\begin{itemize} \item[1)] for $y\in\Omega$, the integral
$$J_{\a}(y)=\int_{\mathbb{R}^{n}}\left|\Delta^{-\a}(x+iy)\right|dx
$$ converges if and only if $\a > \frac{2n}{r} -1.$ In this case,
$J_{\a}(y)=C_{\a}\Delta^{-\a +\frac{n}{r}}(y)$, where $C_{\a}$ is
a constant depending only on $\a$.\item[2)] For $u\in \Omega$, the
integral
$$\int_{\Omega}\Delta^{-\a}(y+u)\Delta^{\nu-\frac nr}(y)dy
$$ converges if and only if $\nu>\frac nr-1$ and $\a> \nu+\frac{n}{r}
-1$, in which case equals $c_\al\Dt^{\nu-\al}(u)$.\item[3)]The
function $F(z)=\Delta^{-\a}(\frac{z+it}{i})$, with $t \in \O$,
belongs to $A_{\nu}^{p}$ if and only if
$$\nu>\frac nr-1\mand \a > \frac 1p(\nu+\frac{2n}{r} -1).$$ In this
case,$$||F||_{A_{\nu}^{p}}=C_{\a,p}\Delta^{-\a +(\nu+\frac
nr)\frac1p}(t).$$ \end{itemize}
\end{lemma}
We refer to the literature for the proof \cite{BBGNPR}. It means
in particular, using (\ref{zero}), that for $p>\tilde p_\nu$ the
function $F(z)=\Delta^{-\frac nr+1}({z+i\be})\in A^p_\nu$ and is
annihilated by $\Box$; so, there is no Hardy inequality for such
values of $p$. In this range of $p$, as mentioned in the
introduction, the Bergman projection $P_\nu$ is not bounded in
$L^p_\nu$, so we have proved easily Theorem \ref{hardyth} for
$p>\tilde p_\nu$. We shall concentrate on the other values of $p$
later on.
\medskip

Let us now recall the following density properties (see
eg \cite{BBPR,GAR}).

\begin{lemma}\label{density} Let $1\leq p <\infty$ and $\nu> \frac
nr -1$. Then, for all $1\leq q\leq \infty$ and $\mu>\frac nr-1$,
the subspace $A^p_\nu \cap A^q_\mu$ is dense in $A^p_\nu$.
Moreover, $A^\infty\cap A^q_\mu $ is dense in $A^\infty$ for the
weak$^*$-$(L^\infty,L^1)$ topology.
\end{lemma}
\Proof
Let us consider the case $p=\infty$, which is the only new part.
If $F\in A^\infty$ the functions $\Delta ^{-\alpha}((\varepsilon
z+i\be)/i)F(z)$ are in $A^p_\mu\,\cap\,A^\infty$ for large values
of $\alpha $, and we clearly have the required property when
$\varepsilon$ tends to $0$ by Lebesgue dominated convergence
theorem. \ProofEnd

\subsection{Integral operators}
For the characterizations of Besov spaces, we shall need some
integral estimates involving Bergman kernel functions. We consider
the following integral operators
\begin{equation}\label{T}
  T_{\nu,\al}F(z)=\Dt^{\al}(\im z)\int_{\cD}B_{\nu+\al}(z,w)F(w)dV_\nu(w),
\end{equation}
 and
 \begin{equation}\label{T+}
T^+_{\nu,\al}F(z)=\Dt^{\al}(\im
z)\int_{\cD}|B_{\nu+\al}(z,w)|F(w)dV_\nu(w),
\end{equation}
 when these integrals make sense.
 Observe that $P_{\nu}=T_{\nu,0}$.

 \begin{lemma}\label{bergtype} Let $\al,\nu,\mu\in\SR$ and $1\le
 p<\infty$.
 Then the following conditions are
 equivalent:
\begin{itemize}
 \item[(a)]
The operator $T^+_{\nu,\al}$ is  well defined and bounded on
$L_{\mu}^{p}(\cD)$. \item[(b)] The parameters satisfy
$\nu+\al>\frac nr-1$ and the  inequalities
\[
\nu p-\mu>(\tfrac nr-1)\max \{1, p-1\}, \quad \al p+\mu\,>(\tfrac
nr-1)\max \{1, p-1\}.\]
 \end{itemize} \end{lemma} \Proof This result is implicit in
 \cite{BT}. For a complete proof see \cite{sehba}.
\ProofEnd
 In particular,
 when $\nu=\mu>\tfrac nr-1$ and when $p>(\mu +\tfrac nr -1)/\mu$, the condition is satisfied for
 $\alpha$ large enough. We remark that, concerning the operators $T_{\nu,\al}$, the
 sufficient conditions for  $L^p_\mu$-boundedness contained in the previous lemma are far from necessary.
 Indeed, we mentioned this in the introduction for the special case of Bergman projections
 (i.e., $\al=0$ and $\mu=\nu$), where other
 methods, which could be generalized to other values of parameters, give additional ranges of boundedness
 (see Remark \ref{R2} below).

 \begin{lemma}\label{bergtypebis} For $\al,\nu\in\SR$, with $\nu>\tfrac nr -1$.
  Then  the operator $T_{\nu,\al}$ (resp. $T_{\nu,\al}^+$) is bounded in
 $L^\infty$ if and only if $\al>\tfrac{n}{r}-1.$\end{lemma}

 \Proof This follows easily from part $3)$
 of Lemma \ref{int-compl} (see details in \cite{sehba}). Remark that now
 we can write $T$ instead of $T^+$, the condition being also necessary for
 $T$.
 \ProofEnd

\subsection{Reproducing formulas}
We will make an extensive use of the following ``\emph{integration by
parts}''. For $\nu>\frac nr-1$, $1\leq p\leq \infty$ and $F\in
A^p_\nu$, $G\in A^{p'}_\nu$, we have the formula
\begin{equation}\label{parts}
\int_{\cD} F(z)\overline G(z) dV_\nu(z)=c_{\nu, m}\int_{\cD}
F(z)\overline {\Box^m G}(z) \Delta^m(\im z) dV_\nu(z).
\end{equation}
Indeed, the formula holds for $p=2$, where it can be obtained
using Plancherel and the Paley-Wiener characterization of $A^2_\nu$ (see eg \cite{FK}).
 The general case follows by density,
using the fact that ${\Box^m G}(x+iy) \Delta^m(y)$ is also in
$L^{p'}_\nu$ by (\ref{mean}). We can now write the following
general reproducing formula. In the next proposition, we write $c$
for some constant that depends on the parameters involved.
\begin{proposition}\label{Repr}
Let $\nu>\frac nr-1$ and $1\leq p\leq \infty$. For all $F\in
A^p_\nu$ we have the formula
\begin{equation}\label{repr}
  \Box^\ell F(z)=c\int_{\cD} B_{\nu+\ell} (z, w)\Box^m
  F(w)\Delta^m(\im w)dV_\nu(w)
\end{equation}
for $m\geq 0$ and $\ell$ large enough so that $B_{\nu+\ell} (z,
\cdot)$ is in $L^{p'}_\nu$. In particular, when $1\leq p< \tilde
p_\nu$, the formula is valid with $\ell=0$.
\end{proposition}
\Proof We can assume that $m=0$. If not, we use (\ref{parts}). It
is true for $p=2$ and $\ell=0$ because of the reproducing property
of the Bergman projection. Derivation under the integral and
\eqref{boxB} gives also the case $\ell>0$. We then use density in
general. \ProofEnd
\begin{corollary} \label{Repr2}Let $1\le p<\tilde p_\nu$ and
$\nu>\frac{n}{r}-1$. Then every $F\in A^p_\nu$ can be written as
\begin{equation}\label{repr2}
   F(z)=\int_{\cD} B_{\nu} (z, w)F(w)dV_\nu (w).\end{equation}
\end{corollary}


We shall state two more results which can be similarly proved by
density and absolute convergence of the involved integrals
(together with Lemma \ref{int-compl} (3) to verify the statements
about the Bergman kernels).

\begin{proposition}
\label{pro1} Let $\nu>\frac nr-1$ and $\al>\frac nr-1$. Then
$B_{\nu+\al}(\cdot, i\be)\in L^1_\nu$, and for all holomorphic $F$
with $\Dt^\al(\im z)F(z)\in L^\infty$ and all $m\geq0$ we have \Be
F(z)\,=\,c\,\int_{\cD} B_{\nu+\al} (z, w)\,\Box^m
F(w)\,\Dt^{\al+m}(\im w)\,dV_\nu(w). \label{Rep4} \Ee
\end{proposition}

\begin{proposition}
\label{pro2} Let $\mu,\nu, \al \in\SR$ and $1\leq p<\infty$
satisfying \[\nu + \al > \tfrac nr-1,\quad
\nu p-\mu>(p-1)(\tfrac nr-1)\mand\mu + \al p > (p-1)(\tfrac
nr-1)-\tfrac nr.\] Then, $\Dt^{\nu-\mu}(\im z)B_{\nu+\al}(z,
i\be)\in L^{p'}_\mu$, and for all holomorphic $F$ with
$\Dt^{\al}(\im z)F(z)\in L^p_\mu$ we have \Be F(z)=\int_{\cD}
B_{\nu+\al} (z, w)\,F(w)\,\Dt^{\al}(\im w)\,dV_\nu(w).
\label{Rep5} \Ee
\end{proposition}

\section{Proofs of Theorems \ref{hardyth} and \ref{duality}}
\setcounter{equation}{0}

\Proofof{Theorem \ref{hardyth}} Let us first assume that $P_\nu$
is bounded, which implies in particular that $p<\tilde p_\nu$,
that is, $B_{\nu} (z, \cdot)$ is in $A^{p'}_\nu$. Then the formula
$$ F(z)=c\int_{\cD} B_{\nu} (z, w)\Box
  F(w)\Delta(\im w)dV_\nu(w)$$
  implies that $F$ is the projection of the function $\Box
  F(w)\Delta(\Im w)\in L^p_\nu$. The Hardy inequality follows from the continuity
  of the operator.
  \medskip

  Next, consider $2<p<\infty$ and assume that the inequality
  (\ref{hardy}) holds. We  can restrict to the
  range $2<p\leq\tilde p_\nu$, since for larger values $p>\tilde
  p_\nu$,
  as we have seen above, the Box operator is
  not injective in $A^p_\nu$, and hence Hardy's inequality does not hold.

Our proof uses Hardy's inequality, not only for the Box operator,
but for its power $\Box ^m$, with $m$ large enough. We shall use
the following lemma.

  \begin{lemma}\label{L_hardy_m}
  Let $\nu>\frac nr-1$ and $2\leq p\leq\tilde p_\nu$.
 Then,
  \Be\label{hardy_mult1}
\bigl\|\Box F\bigr\|_{L^p_{\nu+p}}\,\leq\,C\,\bigl\|\Box^{m+1}
F\bigr\|_{L^p_{\nu+(m+1)p}},\quad\forall\;F\in
A^p_\nu,\;\;\forall\;m\geq 1.
  \Ee
  \end{lemma}

\Proof Using (\ref{parts}) we can write $$ \Box F(z)=c\int_{\cD}
B_{\nu+p} (z,
 w)\Box^{m}
 \left(\Box F(w)\right)\Delta
 ^m(\im w)dV_{\nu+p}(w),$$  since $\Box F\in A^p_{\nu+p}$ and $B_{\nu+p}(\cdot,z)\in
 A^{p'}_{\nu+p}$.
So the inequality (\ref{hardy_mult1}) follows from the fact that
the
 projector $P_{\nu+p}$ is bounded on $L^p_{\nu+p}$ (since
 the condition on $p$ implies $p<\bar p_{\nu+p}$).
\ProofEnd

 So our assumption that Hardy's inequality (\ref{hardy}) holds implies that,
 for all $F\in A^p_\nu$ and all positive integer $m$, we have the inequality
\begin{equation}\label{hardy_mult}
  \int\!\!\int_\cD |F(x+iy)|^p \,\Delta^{\nu- {n \over r}}
(y)\,dx\,dy \, \leq\, C\,\int\!\!\int_\cD \,\bigl|\Delta ^m (y)
\Box^{m}
 F(x+iy)\bigr|^p\,\Delta^{\nu - {n \over r}} (y)\,dx\,dy.
\end{equation}
We want to prove the existence of some constant $C$ such that, for
$f\in L^p_\nu\cap L^2_\nu$, we have the inequality
$$\|P_\nu f\|_{A^p_\nu}\leq C \| f\|_{L^p_\nu}.$$
Consider such an $f$ with $\| f\|_{L^p_\nu}=1$. Call $F:=P_\nu f$.
By Fatou's Lemma, it is sufficient to prove that the functions
$F_\varepsilon (z):=F(z+i \varepsilon \be)$, which belong to
$A^p_\nu$, have norms uniformly bounded. So, using
(\ref{hardy_mult}), it is sufficient to prove that $ \Box^{m}F_\e$
is uniformly in $L^p_{\nu+pm}$ for some $m$, which is a
consequence of the fact that $\Box^{m}F$ itself is in
$L^p_{\nu+pm}$ for some $m$
(see eg \cite[Corol. 3.9]{GAR}). To prove this, we use the identity
$$\Box^{m}F(z)=c\int_{\cD} B_{\nu+m} (z, w)f(w)
 dV_\nu(w),$$
  so that
$\|\Box^mF\|_{L^p_{\nu+pm}}=c\,\|T_{\nu,m}f\|_{L^p_\nu}$, and if
$m$ is sufficient large we conclude from Lemma \ref{bergtype}.
This finishes the proof of Theorem \ref{hardyth}. \ProofEnd
\medskip

 \Proofof{ Theorem \ref{duality}}
 We first consider the case $\tilde p_\nu~'<p<\infty$, for which
 the Bergman kernel $B_\nu (\cdot, w)$ belongs to $A^p_\nu$. So,
 if $F$ is in $A^{p'}_\nu$ and if the associated linear form
 $\Phi(F)$, given by
 $$\langle \Phi(F), G\rangle=\int_\cD G(z) \overline{F(z)}dV_\nu(z)$$
  vanishes on $A^p_\nu$, Corollary \ref{Repr2} implies that $F=0$. Thus,
$A^{p'}_\nu$ is embedded into the dual of $A^{p}_\nu$. Assume that
this embedding is onto, and hence by the closed graph theorem that
it has continuous inverse. Since every $f\in L^{p'}_\nu$ defines
an element of $(A^{p}_\nu)^*$ by
 $G\mapsto \int_\cD G(z) \overline{f(z)}dV_\nu(z)$, by assumption there exists $F\in A^{p'}_\nu$ such that
$$\int_\cD G(z) \overline{f(z)}dV_\nu(z)=\int_\cD G(z) \overline{F(z)}dV_\nu(z),\quad\forall\;G\in A^{p}_\nu$$
with $\|F\|_{A^{p'}_\nu}\leq c\|f\|_{L^{p'}_\nu}$. Taking for $G$
the Bergman kernel, we see that $F$ is the projection $P_\nu f$,
so that $P_\nu f$ maps $L^p_\nu$ continuously into itself.

Conversely, assume that $P_\nu $ is bounded in $L^p_\nu$ (and, by
duality, on $L^{p'}_\nu$). Then we have the
 identity
$$\int_\cD G(z) \overline{f(z)}dV_\nu(z)=\int_\cD G(z) \overline{P_\nu f(z)}dV_\nu(z)$$
for all $f\in L^{p'}_\nu$ and $G\in A^{p}_\nu$. Indeed, use the
fact that this equality is valid in $L^2_\nu$, and density. Since
every functional $\gamma\in (A^p_\nu)^*$ can be expressed by
Hahn-Banach as $G\mapsto \langle G,f\rangle_\nu$ for some $f\in
L^{p'}_\nu$ (with $\|f\|_{L^{p'}_\nu}=\|\gamma\|$), the above
identity shows that the functional can be obtained from $P_\nu
f\in A^{p'}_\nu$. So, under the assumption that $P_\nu $ is
bounded in $L^p_\nu$, the embedding $\Phi:
A^{p'}_\nu\to(A^{p}_\nu)^*$ is an isomorphism.

  It remains to consider the case when $1\leq p\leq \tilde
  p_\nu'$, where we know that the Bergman projection is not
  bounded, and hence we want to show that $\Phi$ is not an isomorphism.
  First, it is easy to see that
  $\Phi$ is not injective when $1\leq p<\tilde p_\nu'$.
  Indeed, in that range we may find a (non-null) function $F\in A^{p'}_\nu$ with $\Box F=0$.
  Now, it follows from \eqref{parts} that
\begin{equation}\label{dualtohardy}
\int_{\cD} G(z)\overline F(z) dV_\nu(z)=c \int_{\cD}
G(z)\,\overline {\Box F(z)}
 \Delta(\im z)dV_{\nu}(z),\quad G\in A^p_\nu,
\end{equation}
which implies $\Phi (F)\equiv0$.

Let us now consider the end-point, $p=\tilde p_\nu'$. If $F$ is in
$A^{p'}_\nu$ then $\Box F$ is in $A^{p'}_{\nu+p'}$ and, by
\eqref{dualtohardy}, the norm of  $\Phi(F)$ is  bounded by the
norm of $\Box F$ in this space. So, if $\Phi$ was an isomorphism,
we would have some constant $C$ independent of $F$  such that
$$\|F\|_{A^{p'}_\nu}\leq C \|\Box F\|_{A^{p'}_{\nu+p'}}.$$
This is exactly Hardy inequality, which is not valid for
$p'=\tilde p_\nu$, concluding the proof of the theorem. \ProofEnd

\medskip

The next corollary, which is implicitly contained in the previous
 proofs, will be used later on.

\begin{corollary} \label {hardy-prec} Let $\nu>\tfrac nr -1$ and $1\leq p<\tpn$,
and assume that the Hardy inequality \eqref{hardy} holds for $(p,
\nu)$. Then, for every positive integer $m$, the mapping
$\Box^m:A^p_\nu\to A^p_{\nu+mp}$ is an isomorphism. In particular,
for all $G\in A^p_{\nu+mp}$ the equation $\Box^m F=G$ has a unique
solution in $A^p_\nu$. Moreover,
$$\|F\|_{A^p_\nu}\leq C \,\|G\|_{A^p_{\nu+mp}},$$
for some constant $C>0$.
\end{corollary}
\Proof  When $2\leq p<\tpn$, by the assumption and Lemma
\ref{L_hardy_m} we have the estimate $\|F\|_{A^p_\nu}\leq
C\,\|\Box^m F\|_{A^p_{\nu+mp}}$, for all $F\in A^p_\nu$, so we
only need to establish the surjectivity of $\Box^m$.
Since by assumption and Theorem \ref{hardyth} the Bergman
projection $P_\nu$ is bounded in $L^p_\nu$, given any $G\in
A^p_{\nu+mp}$, the function $F=P_\nu(\Dt^m(\im\cdot) G)$ belongs
to $A^p_\nu$. Moreover, by the reproducing  formula \eqref{Rep5}
we have \[\Box^m F(z)= \int B_{\nu+m}(z,w)\, G(w)\,\Dt(\im
w)^m\,dV_{\nu}(w)\,=\,c\,G(z),\] which proves the surjectivity.

 For $1\leq p\leq 2$, the
stated result is even simpler; injectivity follows from Proposition \ref{Repr}
 (with $\ell=0$) and surjectivity from the explicit formula
involving the fundamental solution of $\Box$ (see \cite[Prop. 3.1]{B}).
\ProofEnd

\section{Besov spaces of holomorphic functions and duality}
\setcounter{equation}{0}

Throughout this section, given $m\in\SN$, we shall denote
$$\cN_m := \{F \in \cH (\cD) : \Box^m F = 0 \}$$
and set $$\mathcal{H}_m(\cD)=\mathcal{H}(\cD)/ \cN_m.$$ For
simplicity, we use the following notation for the normalized
operator Box operator: we write
\begin{equation}\label{normalized}
\Dt^m\Box^m F(z):=\Dt^m(\im z)\Box^mF(z),\quad z\in\Tu.
\end{equation}
For convenience, we shall use the same notations for holomorphic functions and for
equivalence classes in $\cH_m$. Remark that, for $F\in\cH_m$, we
can speak of the function $\Box^m F$.
Sometimes we shall write $\Box_z^{-m}G$ for the class in $\cH_m$
of all  $F\in\cH(\Tu)$ with $\Box^m F=G$. When
$G\in\cH(\Tu)$ this class is non-empty by the standard theory of PDEs with constant coefficients (see eg \cite{T}).

\subsection{Definition of $\SB^p_\mu(\cD)$}

Given $\mu\in\SR$ and $1\leq p<\infty$ we wish to define a Besov
space $\SB^p_\mu(\cD)$ consisting of holomorphic $F$ so that
$\Dt^m\Box^m F\in L^p_\mu$ for sufficiently large $m$. The
following proposition clarifies the dependence of such spaces on
the parameter $m$.

\begin{proposition}
\label{Bes1} Let $\mu\in\SR$ and $1\leq p<\infty$, and two
integers $0\leq k\leq m$.
\medskip

\noindent (i) If $\Dt^k\Box^kF$ is in $ L^p_\mu$, then
$\Dt^m\Box^mF$ is in $L^p_\mu$ and $\|\Dt^m\Box^mF\|_{L^p_\mu}\leq
C\|\Dt^k\Box^kF\|_{L^p_\mu}$.
\medskip

\noindent (ii) If $\mu+kp>\frac nr-1$ and Hardy's inequality
\eqref{hardy} holds for $(p,\nu=\mu+kp)$, then $\Dt^m\Box^mF\in
L^p_\mu$ implies the existence of $\tF\in\cH(\cD)$ so that
$\Box^m\tF=\Box^m F$ and $\|\Dt^k\Box^k \tF\|_{L^p_\mu}\leq C
\|\Dt^m\Box^m F\|_{L^p_\mu}$.
Moreover the function $\tF$ is uniquely determined modulo $\cN_k$.
\end{proposition}
\Proof Assertion (i) follows from \eqref{mean}. We focus on
assertion (ii). The assumption on Hardy's inequality implies that
$\Box^{m-k}:A^p_{\mu+kp}\to A^p_{\mu+mp}$ is an isomorphism, by
Proposition \ref{hardy-prec}. Thus since  $\Box^mF\in
A^p_{\mu+mp}$, there is a unique $H\in A^p_{\mu+kp}$ with
$\Box^{m-k}H=\Box^mF$. Now we take for $\tF$ any holomorphic
solution of $\Box^k \tF=H$. \ProofEnd

Given $\mu\in\SR$, $1\leq p<\infty$ and $m\in\SN$ we define the
space \[ \SB^{p,(m)}_\mu:=\bigl\{F\in\cH_m(\cD)\mid \Dt^m\Box^m
F\in L^p_\mu\bigr\}
\]
endowed with the norm $\|F\|_{\SB^p_\mu}=\|\Dt^m\Box^m
F\|_{L^p_\mu}$. Observe that each element of $\SB^{p,(m)}_\mu$ is
the equivalence class of all analytic solutions of the equation
$\Box^m F =g$, for some $g\in A_{\mu+mp}^p.$ Thus, the spaces are
null when $\mu+mp\leq \frac nr-1$. By the previous proposition,
when $0\leq k\leq m$ and $\mu+kp>\frac nr-1$, the natural
projection \Be \Ba{lllll} \SB^{p,(k)}_\mu & & \longrightarrow& &
\SB^{p,(m)}_\mu\\
F+\cN_k& & \longmapsto &&F+\cN_m\Ea \label{Bpkm}\Ee is an
isomorphism of Banach spaces, provided Hardy's inequality
\eqref{hardy} holds for the indices $(p,\nu=\mu+pk)$.  This leads
us to the following definition.

\begin{definition}
Given $\mu\in\SR$ and $1\leq p<\infty$, we define
$\SB^{p}_\mu:=\SB^{p,(k_0)}_\mu$ where $k_0=k_0(p,\mu)$ is fixed
by \Be k_0(p,\mu):=\min\bigl\{k\geq0\mid \mu+kp>\tfrac
nr-1\;\mbox{ \small{and Hardy inequality holds for $(p,\mu+pk)$}
}\bigr\}. \label{k0}\Ee
\end{definition}  Observe that $\SB^{p}_\mu=A^{p}_\mu$ if
and only if $k_0(p,\mu)=0$. When $1\leq p\leq2$ we have
$k_0(p,\mu)=\min\{k\geq0\mid \mu+kp>\tfrac nr-1\}$. For $p>2$,
however, the exact value of $k_0(p,\mu)$ depends on Conjecture 1,
and we only have the estimate
\[
k_1(p,\mu)\leq k_0(p,\mu)\leq k_2(p,\mu)
\]
where
\[
\Ba{lll} k_1(p,\mu)& = & \min\bigl\{k\geq0\mid \mu+kp>\tfrac
nr-1\mand p<p_{\mu+kp}\bigr\}\\
k_2(p,\mu)& = & \min\bigl\{k\geq0\mid \mu+kp>\tfrac nr-1\mand
p<\bar p_{\mu+kp}\bigr\} \Ea
\]
A simple arithmetic manipulation shows that $k_1\leq k_2\leq
k_1+1$, and hence $k_0\in\{k_1,k_1+1\}$. Of course, the conjecture
should be $k_0(p,\mu)=k_1(p,\mu)$, and hence we are at most one
unit above the best possible integer in the definition of
$\SB^p_\mu$. Observe also that $k_1(p,\mu)$ and $k_2(p,\mu)$ can
also be written as \Beas k_1 & = & \min\Bigl\{k\geq0\mid
k+\tfrac\mu p>\max\;\bigl\{(\tfrac nr-1)\tfrac1p, \;(\tfrac
nr-1)(1-\tfrac2p)-\tfrac1p, \;(\tfrac
nr-1)(\tfrac12-\tfrac1p)\;\bigr\}\,\Bigr\},\\ & & \\k_2& = &
\min\Bigl\{k\geq0\mid k+\tfrac\mu p>\max\;\bigl\{(\tfrac
nr-1)\tfrac1p, \;(\tfrac nr-1)(1-\tfrac2p)\;\bigr\}\,\Bigr\}.
\Eeas Thus, we have $k_0=k_1=k_2$ when $1\leq p\leq3$.
In the light-cone setting, the improved results about Conjecture 1
mentioned in the introduction imply $k_0=k_1$ for  $1\leq p<3+\e$
for some $\e=\e_{\mu,n}>0$.
\bigskip

In all cases, we can summarize part of the discussion above in the
following proposition.
\begin{proposition}\label{identifies} Let $1\leq p<\infty$, $\mu\in\SR$ and $k\ge k_0(p, \mu)$.
Then \[\Box^k\colon \SB^{p}_\mu \to A^p_ {\mu+k p}\] is an
isomorphism of Banach spaces. In particular, $\SB^p_\mu$ is an
isomorphic copy of $A^p_ {\mu+k_0 p}$, and when $\mu>\tfrac nr-1$
then $\SB^{p}_\mu=A^p_ {\mu}$ for all $1\leq p<\bar{p}_\mu$.
\end{proposition}
\bigskip

Finally we define separately the special family
$$\SB^p:=\SB^p_{-n/r}\,=\,\bigl\{\,F\in\cH(\cD)\mid\Dt^k\Box^k F\in L^p(\cD,d\lambda)\,\bigr\}\,,$$ where
$k$ is sufficiently large and $d\lambda(z)={\Dt^{-\frac{2n}r}(\im
z)}dV(z)$, that is the invariant measure under conformal
transformations of $\cD$. When $n=r=1$, $\SB^p$ is the analog in
the upper half plane of the analytic Besov space studied by
Arazy-Fisher-Peetre, Zhu and others \cite{AFP, AFP1,Z1,R}. These
spaces have also been considered in bounded symmetric domains by
 Yan (for $p=2$), Arazy and Zhu \cite{Y1,A,Z2}.

Some remarkable properties of $\SB^p$, which have been or will be presented
elsewhere, are the following:

\sline(i) $\SB^p\hookrightarrow\SB^q$ when $p\leq q$. This follows
from trivial embeddings of Bergman spaces.

\bline (ii) If  $\frac nr\in\SN$, then  $\SB^p$ is M\"obius
invariant, ie $\|F\circ\Phi\|_{\SB^p}=\|F\|_{\SB^p}$, for all
conformal bijections $\Phi$ of $\cD$, at least when $p>2-\frac
rn$; see \cite{G}. This property fails to be true when $\frac
nr\not\in\SN$, and it is unknown whether it may hold for $1\leq
p\leq 2-\frac rn$ (except in the one dimensional setting;
\cite{AFP}).

\bline (iii) For $b$ analytic in $\cD$, the \emph{small Hankel
operator} $h_b:A^2\to A^2$ is defined by $h_b(f)=P(b{\Ol{f}})$.
If $1\leq p<\infty$, then $h_b$
belongs to the Schatten class $\cS_p$ if and only if $b\in\SB^p$.
See \cite{S, BGNS}.

\medskip

\subsection{Properties of $\SB^p_\mu$: image of the Bergman operator and duality}

 Let $\nu>\frac nr-1$, $1\leq p<\infty$ and $\mu\in\SR$.
When $m$ is large we extend the definition of the Bergman projection
  $P_\nu$ to functions $f\in L^p_\mu$, by letting $P_\nu^{(m)}(f)$
 be the equivalence class (in $\cH_m$) of all
holomorphic solutions of
\[
\Box^mF\,=\,c_{\nu,m}\,\int_{\cD} B_{\nu+m}(\cdot,
w)f(w)\,dV_\nu(w).
\]
The constant $c_{\nu,m}$ is as in \eqref{boxB}, so that if $f\in
L^2_\nu\cap L^p_\mu$ then $P^{(m)}_\nu(f)=P_\nu(f)+\cN_m$, and in
this sense we say that $P_\nu^{(m)}$ is an extension of the
Bergman projection. Observe that $P_\nu^{(m)}$ is well defined and
bounded from $L^p_\mu$ into $\SB_\mu^{p, (m)}$ if and only if
$T_{\nu,m}$ is bounded in $L^p_\mu$, and in particular, by Lemma
\ref{bergtype}, when $p\nu-\mu>\max (1, p-1)(\tfrac nr-1)$ and $m$
is sufficiently large. Moreover, it follows from the reproducing
formulas that the operator is onto. Indeed, by Proposition
\ref{pro2}, every $F\in\SB^{p,(m)}_\mu$ satisfies
\[
\Box^m F(z)\,=\,\int_{\cD} B_{\nu+m}(z, w)\,\Box^m F(w)\,\Dt(\im
w)^m\,dV_\nu(w)
\]
provided $m$ is sufficiently large, from which it follows
$F=cP_\nu^{(m)}(\Dt^m\Box^m F)$. Therefore we have shown the
following result, which partially establishes part 1 of Theorem
\ref{th3}.

\begin{proposition}
\label{surj} Let $\nu>\frac nr-1$, $\mu\in\SR$ and $1\leq
p<\infty$ so that \Be p\nu-\mu>\max (1, p-1)(\tfrac
nr-1).\label{prestr}\Ee If $m$
is sufficiently large (depending only on $p$ and $\mu$) 
then $P_\nu^{(m)}$ maps $L^p_\mu$ boundedly onto
$\SB^{p,(m)}_\mu$.
\end{proposition}

\BR\label{Pmnu} In this proposition it is enough to consider
integers $m$ so that $mp+\mu>\max\{1,p-1\}(\fnr-1)$, since in this
case $T^+_{\nu,m}$ is bounded in $L^p_\mu$ (by Lemma
\ref{bergtype}). The result continues to hold as long as
$T_{\nu,m}$ is bounded in $L^p_\mu$, for which we give a better
range of $m$ and $p$ in Proposition \ref{surjreal} below. Remark that $\Box^k P_\nu^{(m)}=P_\nu^{(m+k)}$. We could as well speak of the projection $P_\nu$ from $L^p_\mu$ onto $\SB^{p}_\mu$. \ER


%

Turning to duality one has the following result.

\begin{proposition} \label{duality-besov} Let $\mu\in\SR$ and $1 < p <
\infty$. For any integers $m_1\geq k_0(p,\mu)$ and
 $m_2\geq k_0(p',\mu)$, the dual space $(\SB_{\mu}^{p})^{*}$
identifies with $\SB_{\mu}^{p'}$ under the integral pairing
\Be\langle F,G\rangle_{\mu,m_1,m_2}=\int_{\cD}
{\Dt^{m_1}\Box^{m_1}}F(z)\overline
{{\Dt^{m_2}\Box^{m_2}}G(z)}
\,dV_\mu(z),\quad F\in \SB_{\mu}^{p},\quad G\in
\SB_{\mu}^{p'}.\label{b2}\Ee
Moreover, modulo a multiplicative constant, the pairing
$\langle\cdot,\cdot\rangle_{\mu,m_1,m_2}$ is independent of $m_1$
and $m_2$ satisfying these inequalities.
\end{proposition}

\Proof The last
statement of the theorem follows from the formula of integration
by parts in \eqref{parts}. Thus, we can assume in \eqref{b2} that
$m_1=m_2=m$, for $m$ as large as desired.

If we denote $\Phi_G(F)=\langle F,G\rangle_{\mu,m,m}$, then it is
clear that $\Phi_G$ defines an element of $(\SB^p_\mu)^*$ and that
the correspondence $G\in\SB^{p'}_\mu\mapsto \Phi_G$ is linear and
bounded. To see the injectivity, consider for each  $w\in\Tu$ the
function $F_w=B_{\mu+m}(\cdot-\bar w)$, which belongs to
$\SB^p_\mu$ if $m$ is sufficiently large (by Lemma
\ref{int-compl}). Then Proposition \ref{pro1} gives, for every
$G\in\SB^{p'}_\mu$, the identity
\[
\Phi_G(F_w)\,=\,c\int_{\cD} B_{\mu+2m}(z-\bar w)\,\overline
{{\Box^{m}}G(z)}\,\Dt^{m}(\im z)\,dV_{\mu+m}(z)\,=\,c\,{\Ol{\Box^m
G(w)}},
\]
(for large $m$), from which the injectivity follows easily.

To see the surjectivity, consider $\ga\in (\SB^p_\mu)^*$. Using
the isomorphism $\Box^m:\SB^p_\mu\to A^p_{\mu+mp}$ (in Proposition
\ref{identifies}) we can define an element $\widetilde\ga\in
(A^p_{\mu+mp})^*$ by $\widetilde\ga(H)=\ga(\Box^{-m}H)$.
The functional $\widetilde\ga$ can be extended to
$(L^p_{\mu+mp})^*$ by Hahn-Banach, and therefore there exists a
function $g\in L^{p'}_\mu$ so that we can write
\[
\widetilde\ga(H)=\int H(z)\,{\Ol{g(z)}}\,dV_{\mu+m}(z),\quad H\in
A^p_{\mu+mp}.
\]
Consequently for every $F\in\SB^p_\mu$
\[
\ga(F)=\widetilde\ga(\Box^mF)=\int 
\Box^m F(z)\,{\Ol{g(z)}}\,dV_{\mu+m}(z).
\]
Next, let $G=P^{(m)}_{\mu+m}(g)$ which for large $m$ defines
element of $\SB^{p'}_\mu$ (by Proposition \ref{surj}). We claim
that $\ga=\Phi_G$. Indeed, when $F\in\SB^p_\mu$\Beas \langle
F,G\rangle_{\mu,m,m}& = &\int
\Box^mF(z)\,{\Ol{\Box^mG(z)}}\,dV_{\mu+2m}(z)\\
& = & c\int \Box^mF(z)\,\left[ \int
B_{\mu+2m}(w,z){\Ol{g(w)}}\,dV_{\mu+m}(w)
\right]\,dV_{\mu+2m}(z)\\
& = &  c\int \left[ \int B_{\mu+2m}(w,z)\,\Box^mF(z)\,\Dt(\im z)^m
\,dV_{\mu+m}(z)\right]{\Ol{g(w)}}\,dV_{\mu+m}(w)\\
\mbox{\small(by Proposition \ref{pro2})}& = &  c\int
\Box^mF(z)\,{\Ol{g(w)}}\,dV_{\mu+m}(w)\,=\,c\,\ga(F),\Eeas where
Fubini's theorem is justified by the boundedness of the operator
$T^+_{\mu+m,m}$ in $L^{p'}_\mu$ when $m$ is sufficiently large.
This establishes the claim, and completes the proof of the
proposition. \ProofEnd

As a special case we obtain the following, which establishes part
2 of Theorem \ref{th3}.

\begin{corollary}\label{ApBp}
Let $\nu>\frac nr-1$ and $1<p\leq2$. Then, $(A_{\nu}^{p})^{*}$
identifies with $\SB_{\nu}^{p'}$ under the integral pairing \Be
\langle F,G\rangle_{\nu,m}=\int_{\cD}F(z)\overline
{\Dt^{m}\Box^{m}G(z)}
\,dV_\nu(z),\quad F\in A_{\nu}^{p},\quad G\in
\SB_{\nu}^{p'},\label{b3}\Ee for any integer $m\geq k_0(p',\nu)$.
\end{corollary}

\Proof Just observe that in this range $k_0(p,\nu)=0$ and
$\SB^p_\nu=A^p_\nu$ (see Proposition \ref{identifies}). \ProofEnd

 \BR \label{dual_largep}We observe that the duality of Bergman spaces is still
open for values of $p$ for which the Hardy inequality is not
valid; that is, we do not know any (non trivial) description of
the spaces $(A^p_\nu)^*$ for $p\geq p_\nu$. \ER

\subsection{The Bloch space $\SB^\infty(\cD)$}

 The definition of analytic Besov space and the properties in
previous sections extend in an analogous way to the case
$p=\infty$, for which $\SB^\infty$ is called \emph{Bloch space}.
In fact, the Bloch space in $\Tu$ was already introduced
 in \cite{bek,Be} and shown to be the dual of $A^1(\cD)$.
Here we recall these results, together with some new facts about
the required number of equivalence classes.

The following inequality is elementary, and can be obtained from
the mean value property of holomorphic functions exactly as in
\cite[Prop. 6.1]{BBPR}, so we omit the proof here.

\begin{lemma}\label{lemma_inf}
Let $\nu\in\SR$. Then \Be \label{mean_inf}
\bigl\|\Dt(\im\cdot)^{\nu+1}\Box
F\bigr\|_{L^\infty}\,\leq\,C\,\bigl\|\Dt(\im\cdot)^\nu
F\bigr\|_{L^\infty},\quad\forall\; F\in\cH(\Tu).\Ee
\end{lemma}



   For every integer $m$ we define a Bloch
type space\[ \SB^{\infty,(m)}\,:=\,\bigl\{F\in\cH_m(\cD)\mid
\Dt^m\Box^m F\in L^\infty\bigr\},
\]
endowed with the norm $\|F\|_{\SB^{\infty,(m)}}=\|\Dt^m\Box^m
F\|_\infty$. We simply write $\SB^\infty(\cD)$ for the space
$\SB^{\infty,(m)}$ with $m=\lceil\fnr-1\rceil$, the smallest
integer greater than $\frac nr-1$.  We have the following
property:

\begin{proposition}
\label{bloch1} For all integers $m\geq k>\frac nr-1$, the natural
inclusion of $\SB^{\infty,(k)}$ into $\SB^{\infty,(m)}$
is an isomorphism of Banach
spaces.
\end{proposition}

\Proof
We may assume $m=k+1$. By Lemma \ref{lemma_inf} \[
\bigl\|\Dt^{k+1}\Box^{k+1}
f\bigr\|_{L^\infty}\,\leq\,C\,\bigl\|\Dt^k
\Box^kf\bigr\|_{L^\infty},\quad f\in\SB^{\infty,(k)}.
\] We want to prove the converse inequality, which is the analogue
of Hardy's inequality for $p=\infty$, that is,
\begin{equation}\label{hardy-infty}
\|\Delta^k\Box^k f\|_\infty\leq C \|\Delta^{k+1}\Box^{k+1}
f\|_\infty,
\end{equation}
for all $k>\frac nr-1 $ and all $f\in\cH(\cD)$ for which the left
hand side is finite. Choosing $\nu>\fnr-1$, we may use Proposition
\ref{pro1} to write
 \begin{equation}\Box^k f\,=\,c\,\int_{\cD} B_{\nu+k}(\cdot-\bar{w})\Box^{k+1}f(w)\,\Dt^{k+1}(\im w)\,dV_\nu(w).
 \label{boxk}
\end{equation} The inequality \eqref{hardy-infty} follows from the
 fact that $\int_{\cD} |B_{\nu+k}(z-\bar{w})|dV_\nu(w)\leq C \Dt^{-k}(\im z)$ by Lemma \ref{int-compl}.

 This implies the injectivity of the mapping. Let us finally prove that the
 mapping is onto. Let $f\in \cH(\cD)$ be such that
 $\Dt^{k+1}\Box^{k+1}f$ is bounded. Then the right hand side of
 \eqref{boxk} defines a holomorphic function, which may be written
 as $\Box^k g$. We prove as before that $\Dt^{k}\Box^{k}g$ is
 bounded. Moreover, $\Box^{k+1}g=\Box^{k+1}f$, which proves the
 surjectivity of the mapping.
\ProofEnd

\BR  Observe that when $k\leq \frac nr-1$ the injectivity of $\SB^{\infty,(k)}\to\SB^{\infty,(m)}$
fails.
Indeed, the function $F(z)=\Dt^{k+1-\frac
nr}(z+i\be)$ belongs to $\SB^{\infty,(k)}$ and is
typically not null in $\cH_k$. However, $F$ is zero in $\cH_m$ for
all $m>k$ since, by \eqref{boxdelta} and \eqref{zero}, we have
$\Box^{k+1} F(z)=c\Box\Dt^{1-\frac nr}(z+i\be)=0.$\ER

\BR \label{rem_bloch} We do not whether for some $k\leq \frac nr-1$ the correspondence  $\SB^{\infty,(k)}\to\SB^{\infty,(m)}$ may be surjective. This question can also be formulated as follows:
\emph{Is it possible that every element $f$ of $\SB^{\infty}$
possesses a representative $g$ such that $\Delta^k\Box^k g$ is
bounded, with $k\leq \frac nr -1$?} This is the analogue of the
question in \eqref{Gconj}, which we shall answer partially in section 4.6. It
seems to us that this problem has not been considered before in the
literature.\ER

We now turn to the boundedness of Bergman operators in $L^\infty$.
As we did in $\S4.2$, when $\nu>\frac nr-1$ we may extend the
definition of the Bergman projection $P_\nu$ to $L^\infty$
functions by letting $P_\nu^{(m)}f$ be the equivalence class (in
$\cH_m$) of all holomorphic solutions of
\[
\Box^m F\,=\,c_{\nu,m}\,\int_{\cD}
B_{\nu+m}(\cdot-\bar{w})f(w)\,dV_\nu(w),
\]
To do this it suffices to consider $m>\fnr-1$, since by Lemma
\ref{bergtypebis} the above integral is always absolutely
convergent and moreover
\[
\|P^{(m)}_\nu
f\|_{\SB^{\infty,(m)}}\,=\,\|T_{\nu,m}f\|_\infty\lesssim
\|f\|_\infty.
\]
Thus $P_\nu^{(m)}$ maps $L^\infty\to\SB^{\infty}$ boundedly. The
mapping is surjective, as every $F\in\SB^\infty$ satisfies (by
Proposition \ref{pro1})
\[
\Box^m F(z)\,=\,\int_{\cD} B_{\nu+m}(z, w)\,\Box^m F(w)\,\Dt(\im
w)^m\,dV_\nu(w)
\]
and therefore, $F=cP_\nu^{(m)}(f)$ with $f=\Dt^m\Box^m F\in
L^\infty$. Hence we have established the following  result.

\begin{proposition}\label{bloch2}
When $\nu, m>\frac nr-1$, the Bergman projection $P_\nu^{(m)}$
maps $L^\infty(\cD)$ continuously onto $\SB^{\infty}$.
\end{proposition}

 \medskip

 Concerning duality, we recall the identification of the Bloch space with the dual of the Bergman
 space $A^1_\nu$.

\begin{theorem}[B\'ekoll\'e, \cite{Be}] \label{A1}Let
$\nu,m > \frac{n}{r}-1$. Then the dual space
 $(A_{\nu}^{1})^*$ identifies with the Bloch space $\SB^\infty$ under the integral pairing
\Be\langle F,G\rangle_{\nu,m}=\int_{\cD}F(z)\overline {\Dt(\im
z)^m\Box^m G(z)}dV_\nu(z),\quad F\in A_{\nu}^{1},\;\;
G\in\SB^\infty\label{b1}.\Ee Moreover, the pairing
$\langle\cdot,\cdot\rangle_{\nu,m}$ is independent of $m>\fnr-1$.
\end{theorem}

The proof is entirely analogous to the one presented in
Proposition \ref{duality-besov}, so we omit it. Let now
$\mu\in\SR$. Since $\Box^m:\SB^1_\mu\to A^1_{\mu+m}$ is an
isomorphism when $\mu+m>\frac nr-1$ (by Proposition
\ref{identifies}), we obtain as a corollary the following duality
statement.
\begin{corollary} Let $\mu\in \SR$ and let $m_1, m_2$ be two integers such that $\mu+m_1
>\frac{n}{r}-1$ and  $m_2
> \frac{n}{r}-1$. Then $(\SB_{\mu}^{1})^*$
identifies with the Bloch space $\SB^\infty$ under the integral
pairing
\[\langle F,G\rangle_{\mu,m_1,m_2}=\int_{T_{\O}}L_{m_1}F(z)\overline
{L_{m_2}G(z)}\Dt^{\mu -\frac{n}{r}}(\Im z)dz,\quad F\in
\SB_{\mu}^{1},\;\; G\in\SB^\infty,\] where $ L_m H(z)=
\Dt^{m}(Imz)\Box_{z}^{m}H(z)$. Again, the pairing
$\langle\cdot,\cdot\rangle_{\mu,m_1,m_2}$ is independent of
 $m_1,m_2$ (modulo a multiplicative constant).\end{corollary}

\subsection{A real analysis characterization of $\SB^p_\mu$}
 We briefly recall the
real variable theory of Besov spaces adapted to the cone that was
developed in \cite{BBGR}.

Following \cite[$\S3$]{BBGR}, we consider a \emph{lattice} $\{\xi_j\}$ in
$\Omega$ and  a sequence $\{\psi_j\}$ of Schwartz functions in
$\mathbb R^n$ such that $\hpsi_j$ is supported in an invariant
ball centered at $\xi_j$ and $\sum_j\hpsi_j=\chi_\O$. In particular, the sets $\supp\hpsi_j$ have the finite intersection property and  the norms $\|\psi_j\|_{L^1(\SR^n)}$ are uniformly bounded. Below we
denote by ${\cS}'_{\partial \Omega}$ the space of tempered
distributions with Fourier transform supported in ${\partial
\Omega}$. Observe that $\Box u=0$ (in $\cS'$) implies $\supp\hat
u\subset\partial \Omega\cup(-\partial\O)$.

\begin {definition}\label{rBes}
Given $\nu \in \SR$ and $1 \leq p< \infty,$ we define
$$B_{\nu}^p\, := \,\bigl\{f \in \cS'(\SR^n) \mid \supp\hf\subset{\Ol\O}\mand ||f||_{B_{\nu}^p}
 < \infty  \bigr\}/\cS'_{\partial \Omega},$$
 where the seminorm is given by
 $$||f||_{B_{\nu}^p}\, := \,\bigl(\,\sum_j\, \Delta^{-\nu} (\xi_j)\,{||f * \psi_j||_p^p}\,\bigr)^{1\over p}.
$$
\end {definition}

\noindent It can be shown that $B^p_\nu$ is a Banach space and the
definition is independent on the choice of $\{\xi_j,\psi_j\}$ (see
\cite[$\S3.2$]{BBGR}). In the 1-dimensional setting $B^p_\nu$
coincides with the classical homogeneous Besov space
$\dot{B}^{-\nu/p}_{p,p}(\SR)$ (of distributions with spectrum in
$[0,\infty)$, modulo polynomials).

%
%
%
%
%

In certain cases one can avoid equivalence classes in Definition
\ref{rBes}, and this will turn into a representation of
$\SB^p_\nu$ as a holomorphic function space. We denote by $\cLL
g(z) = (g, e^{i(z|.)})$, $z\in\Tu$, the \emph{Fourier-Laplace
transform} of a distribution $g$ compactly supported in $\O$
(which defines an analytic function in $\cD$). For convenience, we
write $\Upsilon$ for the set of indices $(p,\nu)$ such that \Be
\nu>-\tfrac nr\mand 1\leq p<\tilde p_\nu,\quad\mbox{or
}\quad\nu=-\tfrac nr \mand p=\tpn=1.\label{tpn1}\Ee Then, in
\cite[Lemmas 3.38 and 3.43]{BBGR} the following result is
shown\footnote{The results in \cite{BBGR} are stated only for
$\nu>0$, but remain valid as long as $(p,\nu)\in\Upsilon$.}.

\begin{lemma}\label{ups}
Let $(p,\nu)\in\Upsilon$. Then if $f\in B^p_\nu$

\bline (i) the series $\sum_j f*\psi_j$ converges in $\cS'(\SR^n)$
to a distribution $f^\sharp$;

\bline (ii) the series $\sum_j\cLL(\widehat f\hpsi_j)(z)$
converges uniformly on compact sets to a holomorphic function in
$\Tu$,  denoted $\cE(f)(z)$, which satisfies
$$\Dt(\im z)^{(\nu+\fnr)/p}\,|\cE(f)(z)|\leq C\,\|f\|_{B^p_\nu},\quad z\in\Tu.$$

\bline In addition, the mappings
\[f\in B^p_\nu\longrightarrow f^\sharp \in \cS'(\SR^n)\mand f\in B^p_\nu\longrightarrow
\cE(f)
\in \cH(\Tu)\] are continuous and
injective, and for every $f\in B^p_\nu$ we have\[
\lim_{{y\to0}\atop{y\in\O}}\cE(f)(\cdot+iy)=f^\sharp\quad\mbox{in
$\cS'(\SR^n)$ and in $\|\cdot\|_{B^p_\nu}$}.
\]
\end{lemma}

From this lemma we can define an isometric copy of $B^p_\nu$ (and
hence of $\SB^p_\nu$) as a holomorphic function space in
$\cH(\cD)$:

\begin{definition}
For $(p,\nu)\in\Upsilon$ we define the holomorphic function space
$$\mathcal{B}_{\nu}^p := \left\{F = \cE f : f \in B_{\nu}^p  \right\},$$
endowed with the norm $||F||_{\mathcal{B}_{\nu}^p} =
||f||_{B_{\nu}^p}.$
\end{definition}

The following properties hold

\Benu

\item[(a)] $\cB^p_\nu=A^p_\nu$ when Hardy's inequality holds for
$(p,\nu)$, and in particular when $\nu>\fnr-1$ and $1\leq p<\bar
p_\nu$ (see \cite[p. 351]{BBGR}).

\item[(b)] $A^p_\nu\hookrightarrow\cB^p_\nu$ when $\nu>\fnr-1$ and
$1\leq p<\tpn$. The inclusion is strict in the 3-dimensional
light-cone when $\nu<1$ and $p_\nu\leq p<\tpn$.

\item[(c)] $\cB^2_0=H^2(\cD)$ (Hardy space). Moreover,
$\left\{\cB^2_\nu=\cLL\Bigl(L^2(\O;\Dt^{-\nu}(\xi)\,d\xi)\Bigr)\right\}_{\nu>-1}$
is the family of spaces introduced by Vergne and Rossi in the
study of irreducible representations of the group of conformal
transformations of $\cD$ (see \cite{VR} or \cite[Ch. XIII]{FK}).

\item[(d)] If $(p,\nu)\in\Upsilon$ then $\Box :
\mathcal{B}_{\nu}^p \rightarrow \mathcal{B}_{\nu+p}^p$ is an
isomorphism of Banach spaces. This is inherited from the
corresponding property in the scale $B^p_\nu$ (see \cite[Th.
1.4]{BBGR}).

\item[(e)] If $(p,\nu)\in\Upsilon$ then $\SB^p_\nu$ can be
identified with $\cB^p_\nu$, in the sense that every
$F\in\SB^p_\nu$ has a (unique) representative $\tF$ in
$\cB^p_\nu$, and moreover $\|F\|_{\SB^p_\nu}\approx
\|\tF\|_{\cB^p_\nu}$. To show this, let $m=k_0(p,\nu)$ so that
$\Box^{m}F\in A^p_{\nu+mp}=\cB^p_{\nu+mp}$ (by (a)). Then use (d)
to define the unique $\tF\in \cB^p_\nu$ such that
$\Box^{m}\tF=\Box^{m}F$.  \Eenu

The assertion in (e) above gives a representation of $\SB^p_\nu$
as a holomorphic function space with no equivalence classes
involved. For example, when $\nu=-n/r$, the space $\SB^1$ can be
represented by the holomorphic function space $\cB^1_{-n/r}$, even in the one-dimensional setting.
\medskip

Using the box operator, this procedure can be easily extended to
all indices $(p,\nu)$ (not necessarily in $\Upsilon$), to
represent $\SB^p_\nu$ with less equivalence classes than
$k_0(p,\nu)$.
 Namely, given
$\nu\in\SR$ and $1\leq p<\infty$, define\Be
k_*=k_*(p,\nu)\,=\,\min\bigl\{k\in\SN\mid
(p,\nu+kp)\in\Upsilon\bigr\}.\label{m*}\Ee Observe that
$k_*(p,\nu)\leq k_0(p,\nu)$, and the inequality is often strict.
In fact, \[ k_*(p,\nu)=\min\bigl\{ k\mid k+\tfrac\nu p>(\tfrac
nr-1)(1-\tfrac2p)-\tfrac 1p\}
\]
(and  $k_*(1,\nu)=\min\{k\mid k+\nu\geq-\tfrac nr\}$).
Then we have the
following result.

\begin{proposition}\label{4.26}
Let $\nu\in\SR$, $1\leq p<\infty$ and $k_*(p,\nu)$ defined as in
\eqref{m*}. Then every $F\in\SB^p_\nu$  has a unique
representative $\tF$, modulo $\cN_{k_*}$, such that $\Box^{k_*}
\tF\in\mathcal{B}^p_{\nu+k_*p}$, and moreover
$\|F\|_{\SB^p_\nu}\approx \|\Box^{k_*}\tF\|_{\cB^p_{\nu+k_*p}}$.
In particular, $\SB^p_\nu$ identifies with the space \Be\{G\in
\cH_{k_*} \mid \Box^{k_*} G\in
\mathcal{B}^p_{\nu+k_*p}\}.\label{cBp2}\Ee
\end{proposition}
\Proof Combine the fact that $\SB^p_{\nu+k_*p}$ identifies with
$\cB^p_{\nu+k_*p}$ (by property (e) above), with the trivial
isomorphism $\Box^{k_*}:\SB^p_\nu\to\SB^p_{\nu+k_*p}$. \ProofEnd



\bigskip

%

We turn now to the identification between the spaces $\SB^p_\nu$
and $B^p_\nu$ via boundary values, as asserted in the
introduction. When $(p,\nu)\in\Upsilon$ the result is immediate
from (e) above.

\begin{corollary}\label{bdryBp}
Let $(p,\nu)\in\Upsilon$. Then

\sline (i) if $F\in \SB^p_\nu$, there exists
$\lim_{{y\to0}\atop{y\in\O}} \tF(\cdot+iy)=f$ in $B^p_\nu$ (and
$\cS'$), for some representative $\tF$ of $F$.

\sline (ii) if $f\in B^p_\nu$, there exists (a unique)
$F\in\SB^p_\nu$ such that $\lim_{{y\to0}\atop{y\in\O}}
F(\cdot+iy)=f$ in $B^p_\nu$.

\sline In either case\[ \tfrac
1c\,\|f\|_{B^p_\nu}\,\leq\,\|F\|_{\SB^p_\nu}\,\leq\,c\,\|f\|_{B^p_\nu}.
\]
\end{corollary}

 The inverse mapping in (ii) is defined by the operator
$f\mapsto F=\cE(f)$. For general parameters $p$ and $\nu$, $\cE f
$ is no longer defined when $f\in B^p_\nu$, but $\cE(\Box^{k_*}f)$
is well-defined and belongs to $\cB^p_{\nu+k_*p}$. Thus, using
Proposition \ref{4.26}, we may consider a new operator
$\mathtt{E}$ from $B^p_\nu$ into $\SB^p_\nu$ by
$$\Box^{k_*}\mathtt E f:= \cE(\Box^{k_*}f).$$
It is easily seen that $\mathtt E: B^p_\nu\to\SB^p_\nu$ is an
isomorphism, which commutes with the Box operator\[
\Box_z^\ell\circ\mathtt E=\mathtt E\circ
\Box^\ell_x,\quad\quad\forall\;\ell\in\SN.\] Moreover, duality can
be expressed through this isomorphism. Recall first that (see
\cite[$\S3.2$]{BBGR})
$$(B^p_\nu)^*=B^{p'}_\nu$$
whenever the definition of the duality pairing is given by
\begin{equation}\label{bracket}
 [ f,
g]_\nu:=\sum_j\langle f, \Box^{-\nu}
g*\psi_j\rangle,\quad\quad f\in B^p_\nu,\quad g\in B^{p'}_\nu.
\end{equation}
On the right hand side the brackets stand for the action of the
distribution $f$ on the conjugate of the given test function, while $\Box^{-\nu}$
is defined on the Fourier side by the multiplication by
$\Delta(\xi)^{-\nu}$. Then, the duality result in Proposition
\ref{duality-besov} can also be obtained from the above
discussion, since when $F=\mathtt  E f\in\SB^p_\mu$, $G=\mathtt E
g\in\SB^{p'}_\mu$ and $m$ is large we have \[ \langle
F,G\rangle_{\mu,m,m}= c_{m,\mu}\,[f,g]_\mu.\]

Finally,
using real variable techniques we are able to improve on the
results in Proposition \ref{surj} concerning the range of $p$ and
number $m$ for which there is boundedness of $P_\nu^{(m)}$ from
$L^p_\mu$ into $\SB^p_\mu$. Below we consider $P_\nu$ as  a
densely defined operator in $L^p_\mu\cap L^2_\nu$.

\begin{proposition}
\label{surjreal} Let $\nu>\frac nr-1$, $\mu\in\SR$ and $1\leq
p<\infty$ so that \Be p\nu-\mu>\max \{p-1,2-p\}\,(\tfrac nr-1).
\label{lower2}\Ee If $k_*=k_*(p,\mu)$ is as in \eqref{m*}, then
$\Box^{k_*}\circ P_\nu$ extends as a bounded surjective
mapping from $L^p_\mu$ onto $\cB^{p}_{\mu+k_*p}$.
\end{proposition}

\BR\label{R2}  As a special case we obtain that, in the range in
\eqref{lower2}, $P_\nu^{(k_0)}$ maps $L^p_\mu$ continuously onto
$\SB^{p}_\mu$, which in particular establishes part 1 of Theorem
\ref{th3}. Equivalently, the operator $T_{\nu,m}$ in $\S2.3$ is
bounded in $L^p_\mu$ for all $m\geq k_0(p,\mu)$; see the
discussion preceeding Proposition \ref{surj}. \ER

When $\mu=-n/r$ the condition \eqref{lower2} produces no
restriction in $p$, and we obtain the following.

\begin{corollary}
For all $\nu>\frac nr-1$ and $1\leq p<\infty$, the operator
$P^{(k_0)}_\nu$ maps $L^p(\cD,d\lambda)$ onto $\SB^p$. Moreover,
$P_\nu$ extends boundedly from $L^1(d\la)$ onto $\cB^1$.
\end{corollary}

\Proofof{Proposition \ref{surjreal}} The continuity follows from a
similar reasoning as in \cite[Prop. 4.28]{BBGR}, where the case
$\mu=\nu$ was proved. For completeness, we sketch here the
modifications of the general case. Given $f\in L^p_\mu\cap
L^2_\nu$, since $P_\nu f\in A^2_\nu$ we can write it, by the
Paley-Wiener theorem, as $P_\nu f=\mathcal{L}g$, for some $g\in
L^2(\O,\Dt^{-\nu}(\xi)d\xi)$. We must show that $\Box^{k_*}P_\nu
f=\mathcal{L}(\Dt^{k_*}g)$ belongs to $\cB^p_{\mu+k_*p}$, or
equivalently that the inverse Fourier transform of the distribution $\Dt^{k_*}g$
belongs to the real space $B^p_{\mu+k_*p}$. Arguing by duality as
in \eqref{bracket}, this is equivalent to show that for all smooth
$\phi$ with compact spectrum in $\O$
\[
\Bigl|\,\langle\Dt^{k_*}g,\Dt^{-\mu-k_*p}\hphi\rangle\,\Bigr|\,\leq\,C\,\|f\|_{L^p_\mu}\,\|\phi\|_{B^{p'}_{\mu+k_*p}}.
\]
By the Paley-Wiener theorem for Bergman spaces (see eg \cite[p.260]{FK})\Beas
LHS& = & \int_\O g(\xi)\,\Dt^{-\mu-k_*(p-1)}(\xi)\,{\overline{\hphi(\xi)}}\;\tfrac{\Dt^\nu(\xi)}{\Dt^\nu(\xi)}\;d\xi\\
& = & \int\!\!\int_\Tu P_\nu f(w)\,{\overline{\cE(\Box^{\nu-\mu-k_*(p-1)}\phi)}}(w)\,dV_\nu(w)\\
\mbox{\small (since $P^*_\nu=P_\nu$) }& = & \langle
f,\cE(\Box^{\nu-\mu-k_*(p-1)}\phi)\rangle_{dV_\nu}\,
\leq\,\|f\|_{L^p_\mu}\,\|\Dt^{\nu-\mu}\cE(\Box^{\nu-\mu-k_*(p-1)}\phi)\|_{L^{p'}_\mu}.
\Eeas If $p>1$ the last norm equals
\[
\|\cE(\Box^{\nu-\mu-k_*(p-1)}\phi)\|_{L^{p'}_{(\nu-\mu)p'+\mu}}.
\]
Under the conditions \eqref{lower2} we have
$A^{p'}_{(\nu-\mu)p'+\mu}=\cB^{p'}_{(\nu-\mu)p'+\mu}$, since
Hardy's inequality holds for the corresponding indices. Thus,
\[
\|\cE(\Box^{\nu-\mu-k_*(p-1)}\phi)\|_{A^{p'}_{(\nu-\mu)p'+\mu}}\,
\approx\,\|\Box^{\nu-\mu-k_*(p-1)}\phi\|_{B^{p'}_{(\nu-\mu)p'+\mu}}\,\lesssim\,
\|\phi\|_{B^{p'}_{\mu+k_*p}},
\]
as we wished to prove. When $p=1$ one must use  instead\[
\|\Dt^{\nu-\mu}\cE(\Box^{\nu-\mu}\phi)\|_{L^\infty}\lesssim
\|\Box^{\nu-\mu}\phi\|_{B^\infty_{\nu-\mu}}\,\approx
\,\|\phi\|_{B^\infty_0}
\]
(see Lemma \ref{binf} below), and conclude again by duality.
The surjectivity of $\Box^{k_*}\circ P_\nu$
follows from the surjectivity of the operator $P^{(m)}_\nu:L^p_\mu
\to\SB^{p,(m)}_\mu$ for large $m$ in Proposition \ref{surj}, since
the spaces $\cB^p_{\mu+k_*p}$ and $\SB^{p,(m)}_\mu$ are related by
isomorphisms. \ProofEnd

\subsection{A real variable characterization of $\SB^\infty$}

For completeness, we give here the real variable characterization of the
Bloch space $\SB^\infty$, starting with the definition of the distribution
spaces $B^\infty_\nu$ introduced in \cite{BBGR}.

\begin {definition}\label{rBloch} For $\nu\in\SR$ we let
$$||f||_{B_\nu^\infty}\, = \,\sup_j\,\Dt(\xi_j)^{-\nu} ||f * \psi_j||_\infty,
\quad f\in S'(\SR^n),$$ and define the space $B_\nu^\infty$ by
$$B^\infty_\nu\, := \,\bigl\{f \in \cS'(\SR^n) \mid \supp\hf\subset{\Ol\O}\mand ||f||_{B_\nu^\infty}
 < \infty  \bigr\}/\cS'_{\partial \Omega}.$$
\end {definition}

\bigskip

The following result is the analogue of Lemma \ref{ups} for
$p=\infty$. The result was not stated in \cite{BBGR}, so we sketch
the proof for completeness.

\begin{lemma}\label{binf}
Let $\nu>\fnr-1$ and $f\in B^\infty_\nu$. Then

\bline (i) $\sum_j f*\psi_j$ converges in $\cS'(\SR^n)$ to a
distribution $f^\sharp$;

\bline (ii) $\sum_j\cLL(\widehat f\hpsi_j)(z)$ converges uniformly
on compact sets of $\Tu$ to a holomorphic function $\cE(f)(z)$,
which satisfies
$$\Dt(\im z)^{\nu}\,|\cE(f)(z)|\leq C\,\|f\|_{B^\infty_\nu},\quad z\in\Tu.$$
\end{lemma}
\Proof By duality, (i) is equivalent to $\cS(\SR^n)\hookrightarrow
B^1_{-\nu}$, which in view of \cite[Prop 3.16]{BBGR} happens if
and only if $\nu>\fnr-1$. Concerning (ii) and reasoning as in the
proof of \cite[Prop 3.43]{BBGR}, it suffices to see that
$\mathcal{F}^{-1}( e^{-( \be|\cdot)}\chi_\Omega )$ belongs to the
space $B^1_{-\nu}$. Using the isomorphism $\Box^{2\nu}$ and the
identity $\cB^1_\nu=A^1_\nu$ this is equivalent to
$\mathcal{L}(\Dt^{2\nu}e^{-( \be|\cdot)}\chi_\Omega
)(z)=c\Dt(z+i\be)^{-2\nu-\fnr}\in A^1_\nu$, which by Lemma
\ref{int-compl} happens if and only if $\nu>\fnr-1$. \ProofEnd

For simplicity we denote $B^\infty=B^\infty_0$, which can be
identified with the Bloch space $\SB^\infty$ as follows.


\begin{proposition}
For all $k>\frac nr-1$, the correspondence
\[
f\in B^\infty \longmapsto \Box^{-k}_z\bigl[\cE(\Box^k f)\bigr]\in
\SB^\infty
\]
is an isomorphism of Banach spaces.
\end{proposition}

\Proof Since $\Box^k f\in B^\infty_k$, by the previous lemma the
function $G:=\cE(\Box^{k}f)$ is holomorphic in $\Tu$ and
$\Delta^k(\im z) G(z)$ is bounded. Thus the equivalence class of
all $F$ such that $\Box^k_z F=G$ belongs to $\SB^\infty$, and the
correspondence $f\mapsto F+\cN_k$ defines a bounded operator from
$B^\infty$ to $\SB^\infty$.

On the other hand, whenever $\nu>\fnr-1$ and $H:=\cE(h)$ is in
$A^1_\nu$, so that $h$ belongs to $B^1_\nu$, one has
$$\int_{\cD}H(z)\overline {\Box^k F(z)}\Delta^k (\Im
z)dV_\nu(z)=[h,f]_\nu.$$ Using the duality identities
$\SB^\infty=(A^1_\nu)^*$ (with the above pairing) and $B^\infty=
(B^1_\nu)^*$ (with the pairing $[\cdot,\cdot]_\nu$), it follows
that the mapping $f\mapsto F$ is an isomorphism, like the mapping
$h\mapsto H$. \ProofEnd

\subsection{Minimum number of equivalence classes: partial results}

Here we consider the question raised in \eqref{Gconj}. We look
first at $p=\infty$ and its equivalent formulation raised in
Remark \ref{rem_bloch}, namely the surjectivity of the mapping
$\SB^{\infty,(k)}\to\SB^{\infty}$ for $k\leq \frac nr-1$. We prove
that it cannot happen at least when $k\leq (\frac nr-1)/2$.

\begin{proposition}\label{NCbloch} Let $k$ be a non negative integer. If, for every $F\in \SB^{\infty}$, there
exists $\tF$ such that $\Delta^{k}\Box^{k} \tF $ is bounded and
$\Box^m\tF=\Box^m F$ for some $m>\frac nr-1$, then necessarily $k>\frac12(\frac
nr -1)$.
\end{proposition}
\Proof Let $m>\frac nr -1$. By the open mapping theorem, if this
property is valid, the natural mapping of $\SB^{\infty,(k)}$ into
$\SB^{\infty,(m)}$, which is surjective, defines an isomorphism from
the quotient space $\SB^{\infty,(k)}/\cN_m$ onto $\SB^{\infty,
(m)}$. So there is some constant $C$ such that, for each $F\in
\SB^{\infty, (m)}$, there exists some $G$ with $\Box^m
G=0$ and
$$\|F+G\|_{\SB^{\infty, (k)}} \leq
C\|F\|_{\SB^{\infty, (m)}}.$$
In particular,
$$|\Box ^k F(x+i\be)+\Box ^k G(x+i\be)|\leq
C\|F\|_{\SB^{\infty, (m)}}.$$ Consider now $F=\cE f$ with $\hf\in C^\infty_c(\O)$, so that
$\|F\|_{\SB^{\infty, (m)}}\leq C\|f\|_{B^{\infty}}$.  Since $\Box ^k F(x+i\be)$ is bounded, the same is valid
for $\Box^k G(x+i\be)$. So we can speak of the Fourier transform of $\Box^k G(x+i\be)$,
whose support is in the boundary of $\Omega$. Let $\phi$ be
a smooth function whose Fourier transform is compactly supported
in $\Omega$, and consider its  scalar product, in the $x$ variable, with the function  $\Box ^k F(x+i\be)+\Box ^k G(x+i\be)$. By the
support condition on  $\hat \phi$ we must have
$\langle \Box ^k G(x+i\be),\phi\rangle=0$. So, the following inequality, valid for all such
$F$, holds
$$\left|\int_{\SR^n}\Box^k F(x+i\be)\overline{\phi(x)} dx\right|\leq C \|f\|_{B^{\infty}}\times \|\phi\|_1.$$
The last inequality can
as well be written as
$$\left|\int_{\SR^n} f(x)\overline{T\phi(x)}dx\right|\leq C \|f\|_{B^{\infty}}\times \|\phi\|_1,$$
where $\widehat{(T\phi)}(\xi)=\Delta(\xi)^k
e^{-(\be|\xi)}\hat \phi(\xi)$. In view of the duality
$(B^1_0)^*=B^\infty$, it is easily seen that this implies the
inequality
\begin{equation}\label{contradiction}
\|T\phi\|_{B^1_0}\leq C \|\phi\|_{1}.
\end{equation}
We want to find a contradiction by choosing specific functions
$\phi$. Assume that $\phi:=\phi_t$ may be written as
$$\phi_t(x)=\sum_{j\in J} r_j(t) a_j e^{i(x|\xi_j)}\eta(x),$$
where $J$ is a finite set of indices, and $\eta$ is a smooth
function whose Fourier transform is supported in a small ball
centered at $0$, in such a way that the functions $\psi_j$ can be
assumed to be equal to $1$ on the support of $\hat \eta (\cdot
-\xi_j)$, for all $j\in J$. Here $r_k(t)$ stands for the
Rademacher function and the parameter $t$ varies in $(0,1)$.
Integrating in $t$ and using Khintchine's Inequality, we have
\Be\label{rademacher}
\int_0^1 \|T\phi_t\|_{B^1} dt\leq C \int_0^1 \|\phi_t\|_{1} dt \leq C'\,\left(\sum_{j\in J}
|a_j|^2\right)^{1/2}\|\eta\|_1.\Ee Let us find a minorant for the
left hand side of \eqref{rademacher}. For every choice of $t$, we
have
$$\|T\phi_t\|_{B^1}= \sum_{j\in J} |a_j|\;
\bigl\|T(e^{i(\cdot|\xi_j)}\eta)\bigr\|_1.$$
Let us take for granted  the existence of some uniform
constants $c_1,c_2>0$ such that
\begin{equation}\label{min-rad}
    \bigl\|T(e^{i(\cdot|\xi_j)}\eta)\bigr\|_1=\Bigl\|\cF^{-1}[\Dt^ke^{-(\be|\cdot)}\hat\eta(\cdot-\xi_j)]\Bigr\|_1\,
    \geq\,\frac1{c_1}\, \Delta(\xi_j)^k\,
e^{-c_2(\be|\xi_j)}\,\|\eta\|_1.
\end{equation}
 Then, \eqref{rademacher} leads to
the existence of  some (different) constant $C$ such that
$$\sum_{j\in J}|a_j|\Delta(\xi_j)^ke^{-c_2(\be|\xi_j)}
\leq C \left(\sum_{j\in J} |a_j|^2\right)^{1/2}.$$ We choose $a_j
=\Delta(\xi_j)^ke^{-c_2(\be|\xi_j)}$ and find  that
$$ \sum_{j\in J}\Delta(\xi_j)^{2k}e^{-2c_2(\be|\xi_j)}\leq C^2$$
uniformly when $J$ varies among finite sets of indices. This allows to have the same estimate for the sum over all indices $j$, that is
 $$ \sum_{j}\Delta(\xi_j)^{2k}e^{-2c_2(\be|\xi_j)}<\infty. $$ By
\cite[Prop. 2.13]{BBGR}  this sum behaves as the
integral
$$\int_\Omega \Delta(\xi)^{2k}e^{-(\be|\xi)}\frac
{d\xi}{\Delta(\xi)^{n/r}},$$ which is finite for $2k>\frac nr-1$.

It remains to prove our claim \eqref{min-rad}, which we do by
using group action as in \cite[(3.47)]{BBGR}. Write $\xi_j=g_j\be$
with $g_j=g^*_j\in G$, and let $\chi_j(\xi)=\chi(g_j^{-1}\xi)$ for
some $\chi\in C^\infty_c(\O)$ with the property that
$\chi_j\equiv1$ in $\supp\hat\eta(\cdot-\xi_j)$, $\forall\;j\in J$
(which we can do by our choice of $\eta$). Consider the function
$\gamma_j$ whose Fourier transform is defined by
$$\widehat{\gamma_j}(\xi):=e^{(\be |\xi)}\Delta(\xi)^{-k}\chi_j(
\xi),$$ so that we can write
$$e^{i(\cdot|\xi_j)}\eta =\gamma_j*T\bigl(e^{i(\cdot|\xi_j)}\eta\bigr),\quad \forall\;j\in J.$$
Thus, it suffices to show that \Be\|\gamma_j\|_1\leq c_1
\Delta(\xi_j)^{-k}e^{c_2(\be|\xi_j)}.\label{dexp}\Ee Now, a change
of variables gives \[ \|\gamma_j\|_1 \, = \,
\bigl\|\,\cF^{-1}\bigl[\,e^{(\be|g_j\xi)}\Dt(g_j\xi)^{-k}\chi(\xi)\,\bigr]\,\bigr\|_1\,
= \,
\Dt(\xi_j)^{-k}\,\bigl\|\,\cF^{-1}\bigl[\,e^{(\xi_j|\cdot)}\Dt^{-k}\chi\,\bigr]\,\bigr\|_1,\]
where in the last equality we have used \eqref{inv-det} and
$g_j^*=g_j$. The $L^1$-norm on the right hand side can be
controlled by a Schwartz norm of $e^{(\xi_j|\cdot)}\Dt^{-k}\chi$,
which leads to \eqref{dexp} using the fact that
$e^{(\xi_j|\xi)}\leq e^{c_2(\xi_j|\be)}$ when $\xi\in\supp\chi$
(see eg \cite[Lemma 2.9]{BBGR}). \ProofEnd

\bigskip

We consider now the same problem for $\SB^p_\nu$, namely the
surjectivity of $\SB^{p,(k)}_\mu\to \SB^{p,(m)}_\mu$ for some
$k<k_0(p,\mu)$. Again, this cannot happen at least if $k$ is
small.

\begin{proposition} \label{nec_incl} Let $\mu\in\SR$ and $k$ be a non negative integer.
If, for every $F\in \SB^{p}_\mu$, there exists $\tF$ such that
$\Delta^{k}\Box^{k} \tF\in L^p_\mu$ and $\Box^m\tF=\Box^m F$ for some
$m\geq k_0(p,\mu)$, then necessarily \Be k+\tfrac \mu
p\,>\,\max\,\Bigl\{(\tfrac nr-1)\tfrac1p\,,\,(\tfrac
nr-1)\,(\tfrac12-\tfrac1p)\,\Bigr\}. \label{condi}\Ee
\end{proposition}

\Proof We must clearly have $\mu+kp>\tfrac nr-1$, since otherwise
$\Box^{k} \tF\in A^p_{\mu+kp}=\{0\}$, which implies $F=0$ (mod
$\cN_m$). We may also assume that
 $k<k_0(p,\mu)$, since otherwise
\eqref{condi} is trivial. In particular, we only need to consider
$p>2$.

The proof is similar to Proposition \ref{NCbloch} with some small
changes. Under the condition in the statement, the inclusion
$\SB^{p,(k)}_\mu/\cN_m\to \SB^{p,(m)}_\mu$ is an isomorphism of
Banach spaces. Hence, for every smooth $f$ with Fourier transform
compactly supported in $\O$, the function $F=\cE(f)$ belongs to $
\SB^{p,(m)}_\mu$ and there exists some $G\in\cH(\cD)$ with $\Box^m
G=0$ so that \Be \|\Dt^k\Box^k(F+G)\|_{L^p_\mu}\lesssim
\|\Dt^m\Box^m F\|_{L^p_\mu}.\label{km}\Ee As before, $\Box^kG$ is
the Fourier-Laplace transform of some distribution supported in
$\bO$. Thus, for all $\hphi\in C^\infty_c(\O)$ we have\Bea
\Bigl|\int_{\SR^n}\Box^k F(x+i\be)\phi(-x)\,dx\Bigr| & = &
\bigl|\Box^k(F+G)(\cdot+i\be)*\phi(0)\bigr|\nonumber\\
&\leq &
\|\phi\|_{p'}\,\bigl\|\Box^k(F+G)(\cdot+i\be)\bigr\|_{L^p(\SR^n)}\label{Fconvphi}.
\Eea Since $\mu+kp>\frac nr-1$ we have
$\bigl\|\Box^k(F+G)(\cdot+i\be)\bigr\|_{L^p(\SR^n)}\lesssim
\bigl\|\Box^k(F+G)\bigr\|_{A^p_{\mu+kp}(\cD)}$ (see e.g.
\cite[Prop. 4.3]{BBGR}). By \eqref{km} and the results in $\S4.4$,
this last quantity is controlled by  \[\|\Box^m
F\|_{A^p_{\mu+mp}}\lesssim \|\Box^m f\|_{B^p_{\mu+mp}}\approx
\|f\|_{B^p_\mu},\] since $m\geq k_0(p,\mu)$. Thus, going back to
\eqref{Fconvphi} we see that
$$\left|\int_{\SR^n} f(x)\overline{T\phi(x)}dx\right|\leq C \|f\|_{B^{p}_\mu}\times \|\phi\|_{p'},$$
where as before
$\widehat{T\phi}(\xi)=\Dt^k(\xi)e^{-(\be|\xi)}\hphi(\xi)$. The
left hand side can be written as a duality bracket $[f, T_\mu
\phi]_\mu $ by letting
$\widehat{T_\mu\phi}(\xi)=\Delta(\xi)^{k+\mu} e^{-(\be|\xi)}
\hphi(\xi)$, and hence we conclude that
\begin{equation}\label{contradiction2}
\|T_\mu\phi\|_{B^{p'}_\mu}\leq C \|\phi\|_{p'}\,.
\end{equation}
  As before, we choose $\phi:=\phi_t$ with
$$\phi_t(x)=\sum_{j\in J} r_j(t) a_j e^{i(x|\xi_j)}\eta(x),$$
where $J$ is a finite set of indices and $\eta$ is a smooth
function with Fourier transform supported in a small ball centered
at $0$ so that  $\psi_j$ can be assumed to be equal to $1$ on the
support of $\hat \eta (\cdot -\xi_j)$, for all $j\in J$.
Integrating in $t$ and using Khintchine's inequality we find that
\begin{equation}\label{rademacher2}
\int_0^1 \|T_\mu\phi_t\|_{B^{p'}_\mu}^{p'}\, dt\leq C \int_0^1
\|\phi_t\|_{p'}^{p'}\, dt\,\leq\,C'\,\left(\sum
|a_j|^2\right)^{p'/2}\|\eta\|_{p'}^{p'}\,,\Ee  while the left hand
side equals
$$ \sum_{j\in J} \Delta(\xi_j)^{-\mu}|a_j|^{p'}
\|T_\mu(e^{i(\cdot|\xi_j)}\eta)\|_{p'}^{p'}.$$ Arguing as in the
proof of \eqref{min-rad} one finds two constants $c_1, c_2$ such
that
$$c_1\|T_\mu(e^{i(\cdot|\xi_j)}\eta)\|_{p'} \geq
\Delta(\xi_j)^{k+\mu} e^{-c_2(\be|\xi_j)}\|\eta\|_{p'}.$$ So,
\eqref{rademacher2} links to the existence of  some  constant $C$
such that
$$\sum_{j\in J}|a_j|^{p'}\Delta(\xi_j)^{kp'+\mu p'-\mu} e^{-c_2(\be|\xi_j)}
\leq C \left(\sum_{j\in J} |a_j|^2\right)^{p'/2}.$$ By the duality
$\ell^r, \ell^{r'}$ with $r=2/p'$ (since we assume $p>2$), we
conclude that
$$ \sum_{j}\Delta(\xi_j)^{r'(kp'+\mu(p'-1))}e^{-c_3(\be|\xi_j)}<\infty,$$
since its partial sums are uniformly bounded. As in the previous
proof, we conclude by a comparison with the corresponding
integral, and find the constraint on parameters in \eqref{condi}.

\ProofEnd

\BR\label{ppsur}
 In the special case $k=0$ we obtain, for $\nu>\fnr-1$ and $m\geq k_0(p,\nu)$, that a necessary condition for the operator $\Box^m:A^p_\nu\to A^p_{\nu+mp}$ to be surjective is
\Be
1\leq p<\,\frac{2(\nu+\frac {n}r-1)}{\frac
nr-1}\,=\,\tpn\,+\,\frac{\nu-1}{\frac nr-1}.
\label{psur}\Ee
When $\nu\leq1$ (in the three dimensional light-cone),  \eqref{psur} is the same necessary condition given in Conjecture \ref{critical}. When $\nu>1$, however, it is
a weaker condition.
\ER

\section{Open Questions}
\setcounter{equation}{0} \setcounter{footnote}{0}
\setcounter{figure}{0}

In this section we pose some questions left open
in this topic, in addition to Conjecture 1. Most questions
concern the spaces $A^p_\nu$ for $p\geq \tpn$, about which we know very little.

\bline (I) \emph{Is the operator $\Box^m:A^{p}_{\nu}\rightarrow
A^p_{\nu+mp}$ onto for some $p\geq\tpn$ and $m\geq k_0(p, \nu)$?}

\sline Equivalently, given a datum $G\in A^p_{\nu+mp}$, does the equation \[
\Box^m F= G\]
have some solution $F$ belonging to the space $A^p_\nu(\Tu)$?

\sline From Remark \ref{ppsur} we only have a negative answer when $p\geq \tpn +(\nu-1)/(\fnr-1)$.

\bline (II) \emph{Is the
operator $\Phi:A^{q'}_{\nu}\rightarrow (A^q_{\nu})^*$ onto for
some $q\leq\tpn'$?}

\bline This question is equivalent to (I) for $p=q'$, using
the duality property  $(A^q_\nu)^*=\SB^{p}_\nu$  in Corollary \ref{ApBp}.

\bline (III) \emph{Is the Box operator injective on
$A_\nu^p$ when $p={\tilde p_\nu}$?}

\bline Injectivity holds when $1\leq p<\tpn$ (by Proposition \ref{Repr}), and fails when $p>\tpn$ (by the explicit example $\Dt(z+i\be)^{-\fnr+1}$). We do not have a conjecture for the endpoint $p=\tpn$.

\bline (IV) \emph{Is the mapping
$\Phi\colon A^{q'}_\nu\to (A^q_\nu)^*$ injective when $q=\tpn'$?}

\bline This is equivalent to (III). In fact, from \eqref{dualtohardy} it easily seen that $ \mbox{Ker }\Phi{|_{A_\nu^{\tpn}}}\,=\,\mbox{Ker}
\Box{|_{A_\nu^{\tpn}}}$.

\bline Our next question stresses further the differences
between the spaces $A^p_\nu$, depending on whether $p<\tpn$ are $p\geq \tpn$:

\bline (V) \emph{Is the space $A^p_\nu$  isomorphic to $\ell^p$ for some $p\geq\tpn$?}

\bline Recall here that the Bergman spaces $A^p_\nu$ are isomorphic to $\ell^p$ in the one dimensional setting. This can be proved as a consequence of the atomic decomposition (see \cite{R}). In \cite{BB2}, atomic decompositions for $A^p_\nu$ are derived when Hardy's inequality holds,
 and they will be developed by the last author in a forthcoming paper also for the spaces $\SB^p_\nu$. We do not know whether
$A^p_\nu$ may be isomorphic to $\ell^p$, or even to $\SB^p_\nu$, when both spaces do not coincide, that is, when Hardy's inequality does not hold.

\bigskip

\bline (VI) \emph{Is it span$\;\{B_\mu(\cdot,w)\mid w\in\Tu\}$ dense in $A^p_\nu(\Tu)$ for $p\geq \tpn$ and $\mu$ sufficiently large?}

\bline The validity of this result was wrongly stated in \cite[Corollary 5.4]{BBPR} in the light-cone setting. As we show below (see also \cite[Lemma 5.1]{BBPR}), the density holds when the projection $P_\mu$ is bounded in $L^p_\nu$, but this restricts
$p$ to be smaller than $\tpn$ (since $P_\mu^*=T_{\nu,\mu-\nu}$ must also be bounded in $L^{p'}_\nu$).


\begin{proposition}\label{density} Let $\nu>\frac nr-1$.
Assume that $p$ and $\mu$ are so that
$P_\mu$ extends as a bounded operator in $L^p_\nu$.
Then $A^p_\nu$ is the
closed linear span of the set $\{B_\mu (., w), w\in \cD\}.$
\end{proposition}

\Proof
The boundedness of $P_\mu$ in $L^p_\nu$ already implies that $B_\mu(\cdot,i\be)\in A^p_\nu$.
We take for granted the fact that
$P^*_\mu=T_{\nu,\,\mu-\nu}$ (with respect to $\langle\cdot,\cdot\rangle_{dV_\nu}$).
To establish the proposition it suffices to prove that, for
$f\in L^{p'}_\nu$ such that
\begin{equation}\label{density11}
\langle f, B_\mu(\cdot,w)\rangle_\nu=0,\quad\forall\;w\in\Tu,
\end{equation}
we have also $\langle f, F\rangle_\nu=0$ for all $F$ in a dense subset of
$A^p_\nu$. Now \eqref{density11} is the same as $T_{\nu,\,\mu-\nu}(f)(w)=0$, by definition of this operator. Thus, if $F\in A^p_\nu\cap A^2_\mu$, using the claim above we have
\[
\langle f, F\rangle_\nu= \langle f, P_\mu F\rangle_\nu=\langle P^*_\mu(f), F\rangle_\nu=0.\]
Finally, we establish the claim, that is $P^*_\mu=T_{\nu,\,\mu-\nu}$.
For $f,g\in C_c(\Tu)$ we have to justify the exchange
of order of integration in
\Beas
\langle P_\mu(g),f\rangle_\nu& = &
\int_\cD \Bigl[\int_\cD B_\mu (z,
w)g(w)dV_\mu (w)\Bigr]\,{\Ol{f(z)}}\,dV_\nu(z)\\
& = & \int_\cD g(w)\,\Bigl[\int_\cD {\Ol{B_\mu (w,
z)f(z)}}dV_\nu (z)\Bigr]\,dV_\mu(w)\,=\,\langle g,T_{\nu,\,\mu-\nu}f\rangle_\nu. \Eeas
but this follows from $$\int_\cD \int_\cD  |B_\mu (z, w)||g(w)|\,dV_\mu(w)|f(z)|\,dV_\nu (z)\leq \bigl\|T^+_{\mu,0}|g|\bigr\|_{L^2_\mu}\,
\bigl\|\Dt^{\nu-\mu}|f|\bigr\|_{L^2_\mu}<\infty,$$
using the fact that the operator $T^+_{\mu,0}$ with kernel $|B_\mu (z, w)|$ is bounded on $L^2_\mu$.
\ProofEnd

\bibliographystyle{plain}

\end{document}